\documentclass[11pt]{amsart} 
\usepackage{fullpage, amsmath, amssymb, young}
\usepackage[all]{xypic}

\title{Poisson traces for symmetric powers of symplectic varieties}
\author{Pavel Etingof and Travis Schedler}
\usepackage{url}            
\usepackage{graphicx}       

\numberwithin{equation}{subsection}
\theoremstyle{definition}
\newtheorem{theorem}[equation]{Theorem}
\newtheorem{lemma}[equation]{Lemma}

\newtheorem{corollary}[equation]{Corollary}

\newtheorem{example}[equation]{Example}

\newtheorem{conjecture}[equation]{Conjecture}
\newtheorem{remark}[equation]{Remark}

\newtheorem{claim}[equation]{Claim}

\newcommand\Tline{\rule{0pt}{2.6ex}}

\newcommand{\Span}{\operatorname{Span}}
\newcommand{\Id}{\operatorname{Id}}

\newcommand{\ad}{\operatorname{ad}}

\newcommand{\HH}{\mathsf{HH}}
\newcommand{\Hilb}{\operatorname{Hilb}}
\newcommand{\HHH}{\mathcal{H}}
\newcommand{\Stab}{\operatorname{Stab}}
\newcommand{\symm}{\text{symm}}

\newcommand{\pt}{\text{pt}}
\newcommand{\HP}{\mathsf{HP}}

\newcommand{\ZZ}{\mathbb{Z}}

\newcommand{\iso}{{\;\stackrel{_\sim}{\to}\;}}
\newcommand{\C}{\mathbb{C}}

\newcommand{\cO}{\mathcal{O}}
\newcommand{\caD}{\mathcal{D}}

\newcommand{\hh}{\mathfrak{h}}

\newcommand{\Hom}{\text{Hom}}
\newcommand{\Aut}{\operatorname{Aut}}

\newcommand{\Ext}{\text{Ext}}

\newcommand{\SL}{\mathsf{SL}}
\newcommand{\GL}{\mathsf{GL}}
\newcommand{\Sp}{\mathsf{Sp}}

\newcommand{\gr}{\operatorname{\mathsf{gr}}}
\newcommand{\Spec}{\operatorname{\mathsf{Spec}}}
\newcommand{\Spf}{\operatorname{\mathsf{Spf}}}
\newcommand{\Weyl}{\mathsf{Weyl}}
\newcommand{\onto}{\twoheadrightarrow}
\newcommand{\into}{\hookrightarrow}
\newcommand{\Sym}{\operatorname{\mathsf{Sym}}}

\newcommand{\CC}{{\Bbb C}}
\newcommand{\ev}{\operatorname{ev}}
\newcommand{\dblparens}[1]{(\!(#1)\!)}
\newcommand{\dblsqbrs}[1]{[\![#1]\!]}
\begin{document}
\date{2011}
\begin{abstract}
  We compute the space of Poisson traces on symmetric powers of affine
  symplectic varieties. In the case of symplectic vector spaces, we
  also consider the quotient by the diagonal translation action, which
  includes the quotient singularities $T^*\CC^{n-1}/S_n$ associated to
  the type $A$ Weyl group $S_n$ and its reflection representation
  $\CC^{n-1}$.  We also compute the full structure of the natural
  $\caD$-module, previously defined by the authors, whose solution
  space over algebraic distributions identifies with the space of
  Poisson traces.  As a consequence, we deduce bounds on the numbers
  of finite-dimensional irreducible representations and prime ideals
  of quantizations of these varieties.  Finally, motivated by these
  results, we pose conjectures on symplectic resolutions, and give
  related examples of the natural $\caD$-module.  In an appendix, the
  second author computes the Poisson traces and associated
  $\caD$-module for the quotients $T^*\CC^{n}/D_n$ associated to type
  $D$ Weyl groups. In a second appendix, the same author provides a
  direct proof of one of the main theorems.
\end{abstract}
\maketitle
\tableofcontents
\section{Introduction}
\subsection{Poisson traces for $S^n Y$ and type $A$ Weyl group
  quotient singularities}\label{ss:typea-int}
Given a Poisson algebra $A$ over $\CC$, a \emph{Poisson trace} is a
functional $A \rightarrow \CC$ which annihilates $\{A, A\}$. These may
also be viewed as functionals on the associated Poisson variety $\Spec
A$ which are invariant under Hamiltonian flows. The space of such
traces is the dual, $\HP_0(A)^*$, to the \emph{zeroth Poisson
  homology}, $\HP_0(A) := A /\{A, A\}$, which also coincides with the
zeroth Lie homology.

Given an affine variety $Y = \Spec A$, we use the notation $\cO_Y :=
A$.  Let $S^n Y := Y^n / S_n = \Spec \Sym^n A$ be the $n$-th symmetric
power of $Y$.  By a symplectic variety, we always mean a smooth
symplectic variety. Let the symbol $\&$ denote the tensor product in
the symmetric algebra.

Our first main result is:
\begin{theorem}\label{symphp0thm}
  Let $Y$ be an affine symplectic variety.  Then, there is a canonical
  isomorphism of graded algebras,
\begin{gather} \label{e:symphp0thm}
  \Sym(\HP_0(\cO_Y)^*[t]) \iso \bigoplus_{n \geq 0} \HP_0(\cO_{S^n
    Y})^*, \\ \notag \phi \cdot t^{m-1} \mapsto \Bigl( (f_1 \&
  \cdots \& f_m) \mapsto \phi(f_1 \cdots f_m) \Bigr),
\end{gather}
wherein the grading is given by $|\HP_0(\cO_{S^n Y})^*| = n$ (on both
sides of the isomorphism), and $|t| = 1$.
\end{theorem}
This will be proved in \S \ref{typeas}.  It is well known that, if $Y$
is connected, then $\HP_0(\cO_Y) \cong H^{\dim Y}(Y)$, the top
cohomology of $Y$, via the isomorphism $[f] \mapsto f \cdot
\operatorname{vol}_Y$, where $\operatorname{vol}_Y$ is the canonical
volume form (i.e., the $\frac{1}{2} \dim Y$-th exterior power of the
symplectic form). We can write the above more explicitly using the
coefficients $a_n(i)$ which give the number of $i$-multipartitions
of $n$ (i.e., collections of $i$ ordered partitions whose sum of sizes is $n$), i.e.,
\begin{equation}
\prod_{m \geq 1} \frac{1}{(1-t^m)^i} = \sum_{n \geq 0} a_n(i) \cdot t^n.
\end{equation}
%%NOTE: as above, ``connected'' here and below are because, otherwise, one
%%could have connected components with dimension less than the dimension of
%%the overall variety. Then the degree ``dim Y'' appearing below would have
%%to be modified.
\begin{corollary}\label{c:mpfla} If $Y$ is connected, then 
$\dim \HP_0(\cO_{S^n Y}) = a_n(\dim H^{\dim Y}(Y))$.
\end{corollary}

We may derive a relationship with Hochschild homology of a
quantization, as follows.  Given an associative algebra $B$, recall
that the zeroth Hochschild homology is $\HH_0(B) := B / [B,B]$.
Recall also that a (formal) deformation quantization of a Poisson algebra
$\cO_X$ is an associative algebra $(A_\hbar, \star)$ over
$\CC\dblsqbrs{\hbar}$ which is isomorphic to $\cO_X\dblsqbrs{\hbar} := \{\sum_{m
  \geq 0} f_m \cdot \hbar^m \mid f_m \in \cO_X\}$ as a
$\CC\dblsqbrs{\hbar}$-module, and such that, for $a, b\in \cO_X$, the
deformed multiplication satisfies $a \star b = ab + O(\hbar)$ and $a
\star b - b \star a = \hbar \{a,b\} + O(\hbar^2)$.  Here $f = g +
O(\hbar^i)$ means that $f-g \in \hbar^i A_\hbar$.
 Note that $A_\hbar$ and $A_\hbar[\hbar^{-1}] \supset
A_\hbar$ are filtered by powers of $\hbar$: $F^i A_\hbar[\hbar^{-1}] =
\hbar^i A_\hbar$. Hence, $\HH_0(A_\hbar[\hbar^{-1}])$ inherits the
filtration.

\begin{corollary}\label{defcor1} Let $Y$ be an affine symplectic variety, 
  and $A_\hbar$ be any deformation quantization of $\cO_{Y}$ (so
  $\Sym^n A_\hbar$ is a deformation quantization of $\cO_{S^n
    Y}$). Then, the natural surjection $\HP_0(\cO_{S^n Y})\dblparens{\hbar}
  \onto \gr \HH_0(\Sym^n A_\hbar[\hbar^{-1}])$ is an
  isomorphism. 
\end{corollary}
The above corollary uses the computation of $\HH_0(\Sym^n
A_\hbar[\hbar^{-1}])$ from \cite{EO}.  For arbitrary quantizations of
$\Sym^n \cO_Y$ (not necessarily of the form $\Sym^n A_\hbar$), we can
still deduce
\begin{corollary}\label{defcor2}
  Let $Y$ be a connected affine symplectic variety and $B_\hbar$ be an
  arbitrary deformation quantization of $\cO_{S^n Y}$.  Then,
  $\dim_{\CC\dblparens{\hbar}} \HH_0(B_\hbar[\hbar^{-1}]) \leq a_n(\dim
  H^{\dim Y}(Y))$.
\end{corollary}
\begin{example}
  In the case that $Y$ is a connected smooth surface and $H^1(Y) = 0$,
  if $A_\hbar$ is any quantization of $\cO_Y$, then \cite[\S 6]{EO}
  constructs the universal formal deformation of $A_\hbar[\hbar^{-1}]^{\otimes
    n} \rtimes S_n$, called $A_\hbar[\hbar^{-1}](n,c,k)$, over
  $\CC\dblparens{\hbar}\dblsqbrs{H^2(Y) \oplus \CC}$, i.e., $\CC\dblparens{\hbar}$-valued
  functions on the formal neighborhood of zero in $H^2(Y) \oplus
  \CC$. Here, $c := (c_1, \ldots, c_{\dim H^2(Y)})$ and $k$ denote
  bases for linear functions on $H^2(Y)$ and $\CC$, respectively, so we
  write $\CC\dblparens{\hbar}\dblsqbrs{H^2(Y) \oplus \CC} = \CC\dblparens{\hbar}\dblsqbrs{c,k}$.  The
  deformation $A_\hbar[\hbar^{-1}](n,c,k)$ is topologically free over $\CC\dblparens{\hbar}\dblsqbrs{c,k}$, i.e., isomorphic to $A\dblparens{\hbar}\dblsqbrs{c,k}$ as a $\CC(\hbar))\dblsqbrs{c,k}$-module.

  Suppose in addition that there exists an integral form of this
  algebra, i.e., a $\CC\dblsqbrs{\hbar}\dblsqbrs{c,k}$-subalgebra
 $A_\hbar(n,c,k) \subseteq A_\hbar[\hbar^{-1}](n,c,k)$
 satisfying
\begin{enumerate}
\item[(a)] $A_\hbar(n,c,k)$ is topologically free over
 $\CC\dblsqbrs{\hbar}\dblsqbrs{c,k}$, and
\item[(b)] $A_\hbar(n,c,k)[\hbar^{-1}] \cong A_\hbar[\hbar^{-1}](n,c,k)$ as algebras over $\CC\dblparens{\hbar}\dblsqbrs{c,k}$.
%\item[(c)] $A_\hbar(n,c,k)/(\hbar) \cong \cO_{S^n Y}[[H^2(Y) \oplus \CC]]$ as
%algebras over $\CC[[H^2(Y) \oplus \CC]]$.
\end{enumerate}
 Such an $A_\hbar(n,c,k)$ exists in all
  examples we know. Then, the subalgebra $e A_\hbar(n,\hbar c,\hbar k) e$, where
  $e = \frac{1}{n!} \sum_{\sigma \in S_n} \sigma \in \CC[S_n]$ is the
  symmetrizer, is a quantization of $\cO_{S^n Y}\dblsqbrs{c,k}$, and the above corollary applies to show that $\HH_0(e
  A_\hbar(n,\hbar c,\hbar k)[\hbar^{-1}] e)$ is generated by at most $a_n(\dim
  H^2(Y))$ elements. In particular, if one specializes at any values
  of $c$ and $k$, then one obtains a deformation quantization of $A$
  and the dimension of the resulting zeroth Hochschild homology as
  a vector space over $\CC\dblparens{\hbar}$ is at most $a_n(\dim H^2(Y))$.

  Note that this is essentially a global version of the Cherednik
  algebra associated to $S_n$: when one replaces $Y$ by $\CC^2$, one
  can recover the Cherednik algebra associated to the Weyl group $S_n$
  from the above (more precisely, one recovers the usual Cherednik algebra
  tensored by $A_\hbar$, since the Cherednik algebra itself involves
  deforming $A_\hbar^{\otimes n-1} \rtimes S_n$, corresponding
  to the reflection representation of $S_n$).  One can conjecture
  that, parallel to Corollary \ref{defcor6}, in fact $\HH_*(e
  A_\hbar[\hbar^{-1}](n,c,k) e) \cong \HH_*(\Sym^n
  A_\hbar[\hbar^{-1}])\dblsqbrs{c,k} \cong H^{2n-*} (\Hilb^n
  Y)\dblparens{\hbar}\dblsqbrs{c,k}$, where $\Hilb^n Y$ is the Hilbert scheme of the
  surface $Y$; see also \S \ref{ss:csr} below. The second isomorphism
  here follows by comparing \cite[Corollary 3.3]{EO} with the known
  cohomology of the Hilbert scheme, because $A_\hbar[\hbar^{-1}]$ is
  an infinite-dimensional simple algebra with trivial center.
\end{example}

% Since distinct finite dimensional irreducible representations of $B$
% are linearly independent, we deduce also
% \begin{corollary}\label{repcor} Assume $X$ is affine and
%   symplectic. Let $B$ be an arbitrary quantization of $\cO_{S^n X}$.
%   Then, $B$ has at most $a_n(\dim H^{\dim X}(X))$ nonisomorphic finite
%   dimensional irreducible representations.
% \end{corollary}
Given an affine Poisson variety $X$ such that $\cO_X$ is nonnegatively
graded and equipped with a Poisson bracket of degree $-d$, one defines
a filtered quantization to be a filtered algebra $B$ over $\CC$ such
that $\gr B = \cO_X$, $[B_{\leq i}, B_{\leq j}] \subseteq B_{\leq
  i+j-d}$, and for $a \in B_{\leq i}$ and $b \in B_{\leq j}$, $\{\gr_i
a, \gr_j b\} = \gr_{i+j-d} [a, b]$.

In the case $X = V$ is a symplectic vector space, the standard
quantization is given as follows: Write $V = U \oplus U^*$ where $U$
and $U^*$ are complementary Lagrangians.  Let $x_1, \ldots, x_n$ be a
basis of $U^* \subseteq \cO_U$ and $\partial_1,
\ldots, \partial_n \in U$ be the dual basis.  Then, the standard
quantization is the ring of differential operators $\caD_U$ filtered
by the Bernstein filtration,\footnote{Since our groups $G$ will be of
  the form $G<\GL(U) < \Sp(V)$, one could alternatively use the order
  filtration and change the grading on $\cO_V$ accordingly; this would
  have the effect of halving the degrees appearing in $\HP_0$ and
  $\HH_0$.} where $(\caD_U)_{\leq k}$ is spanned by elements of the
form $x_{i_1} \cdots x_{i_j}
\partial_{i_{j+1}} \cdots \partial_{i_{\ell}}$, for $\ell \leq k$.  In
other words, this is the filtration generated by
$|x_i|=|\frac{\partial}{\partial x_i}| = 1$. Then, $\gr \caD_U \cong
\cO_V$ (with $d=2$).  (The algebra $\caD_U$ is also known as the Weyl
algebra of $V$.)  Similarly, one could consider the deformation
quantization $\caD_{U,\hbar}$, which as a $\CC\dblsqbrs{\hbar}$-module is
isomorphic to $\caD_{U}\dblsqbrs{\hbar}$, but with the commutation relations
multiplied by $\hbar$: namely, $\caD_{U,\hbar}$ is generated by $x_i$ and $p_i$ with relations $[p_i, x_j] = \hbar \delta_{ij}$ (one can think of $p_i$ as $\hbar \frac{\partial}{\partial x_i}$).
In the case $X = V/G$, then $\caD_{U}^G$ and $\caD_{U,\hbar}^G$ are
filtered and deformation quantizations of $\cO_V^G = \cO_{V/G}$.

Similarly to the preceding theorem, we can consider quantizations of
 $X = S^{n+1} V$, for $V$ a symplectic vector space.  In this case,
we have a decomposition $S^{n+1} V = V^n/S_{n+1} \times V$, where the second
factor is the diagonally embedded $V$, and the $S_{n+1}$ action on
$V^n$ is by the identification $V^n \cong (\CC^n \otimes V)$, where $\CC^n$ is the reflection representation and $V$ is a trivial representation.
So,
$\HP_0(\cO_{S^{n+1} V}) = 0$, since $\HP_0(\cO_V) = 0$.  On the other hand:
\begin{theorem}\label{snp1thm}
  $\HP_0(\cO_{V^n / S_{n+1}})^* \cong \CC$, spanned by the
  augmentation map $\cO_{V^n} \to \CC$.
\end{theorem}
As we will see, Theorem \ref{symphp0thm} reduces, in a sense, to the
above theorem, using $\caD$-modules to localize the problem.  An
elementary proof of the above theorem, that does not require anything
in the main body of the paper, is provided in Appendix
\ref{s:direct-proof-snp1thm}.
\begin{remark} \label{r:ReSc} In \cite{ReScmat}, building on seminal
  work of Mathieu \cite{Matso}, the second author computes more
  generally some of the structure of $\HP_0(\cO_{V^n/S_{n+1}},
  \cO_{V^n}) := \cO_{V^n} / \{\cO_{V^n}^{S_{n+1}}, \cO_{V^n}\}$.  This
  is an $S_{n+1}$-representation whose invariants are
  $\HP_0(\cO_{V^n/S_{n+1}})$.  In \cite{ReScmat}, the argument used
  here is generalized to show, among other things, that the isotypic
  part of $\HP_0(\cO_{V^n/S_{n+1}}, \cO_{V^n})$ corresponding to Young
  diagrams with at most $\dim V + 1$ boxes below the top row coincide
  with the same isotypic part of the subspace of the free Poisson
  algebra on $n$ variables $z_1, \ldots, z_n$ which has degree one in
  each variable, with $S_{n+1}$ action by the reflection
  representation, and with grading given by twice the number of pairs
  of brackets $\{-,-\}$ which appear. Then, Theorem \ref{snp1thm}
  above follows from the fact that the only multilinear Poisson
  polynomial in $z_1, \ldots, z_n$ which is symmetric in all the
  variables is the product $z_1 \cdots z_n$.  The result of
  \emph{op.~cit.}~also implies that the reflection representation
  $\mathfrak{h}$ of $S_{n+1}$ does not occur, and that the isotypic
  component of $\wedge^2 \mathfrak{h}$ occurs with multiplicity
  $\lfloor \frac{n}{2} \rfloor$.  In terms of affine symplectic
  varieties $Y$, these results translate into information about the
  structure of $\HP_0(\cO_{\Sym^n Y}, \cO_{Y^n})$ as an
  $S_n$-representation; see Remark \ref{r:Mphi}.
\end{remark}
In fact, in \S \ref{typeas}, we will deduce the above two theorems
from a more general result (Theorem \ref{sympdmthm}) on the
$\caD$-module $M(X)$ from \cite{ESdm}, which is essentially the
quotient of $\caD_{S^n Y}$ by the right ideal generated by Hamiltonian
vector fields.  See \S \ref{ss:typea-int-dmod} for the statements.
% Let $\caD(X)$ denote the algebra of differential operators on $X$,
% filtered by order.  In the case that $X$ is a Lagrangian subspace of a
% symplectic vector space $U$ and $G < Sp(U)$, \cite{AFLS} computed the
% space $\HH_0(D(X)^G)$.  Using these results, we deduce the
\begin{corollary}\label{defcor3}
  Let $U \subseteq V$ be a Lagrangian subspace.  Then, the natural
  surjection $\HP_0(\cO_{V^n}^{S_{n+1}}) \onto \gr
  \HH_0(\caD_{U^n}^{S_{n+1}})$ is an isomorphism, and both are
  isomorphic to $\CC$.
\end{corollary}
% [[Are all quantizations of $\cO_{V^n}^{S_{n+1}}$ are necessarily
% isomorphic?]]
\begin{remark}
%We may also deduce the following result in noncommutative algebra:
%\begin{corollary}\label{defcor4}
  As a consequence, if $B$ is any filtered quantization of
  $\cO_{V^n}^{S_{n+1}}$, then, $\dim \HH_0(B) \leq 1$, and hence $B$
  admits at most one finite-dimensional irreducible
  representation. However, when $\dim V > 2$, we do not know if
  filtered quantizations not isomorphic to $\Weyl(V^n)^{S_{n+1}}$
  exist, and for the latter the zeroth Hochschild homology was already
  computed in \cite{AFLS} (we discuss the case $\dim V = 2$ below).
%\end{corollary}
\end{remark}
In the case $V = \CC^2$, then $V^n = \mathfrak{h} \oplus
\mathfrak{h}^*$, where $\mathfrak{h} \cong \CC^n$ is the reflection
representation of $S_{n+1}$, viewed as a type $A_n$ Weyl group.  In
this case, the theorem specializes to
\begin{corollary}\label{defcor5}
  Let $\mathfrak{h} \cong \CC^{n}$ be the reflection representation of
  the Weyl group $S_{n+1}$ of type $A_{n}$. Then
\begin{equation}
\HP_0(\cO_{(\mathfrak{h} \oplus \mathfrak{h}^*)/S_{n+1}}) \cong \CC \cong
\HH_0(\caD(\mathfrak{h})^{S_{n+1}}).
\end{equation}
\end{corollary}
This corollary was verified by computer by Justin Sinz for small
values of $n$; the cases $n \leq 2$ and $n = 3$ are also proved in
\cite{AlFo} and \cite{Bu}.

More generally, we can extend the corollary to the case of spherical
rational Cherednik algebras associated to $S_{n+1}$. Recall (see,
e.g., \cite{EGsra}) that these are certain filtered algebras $B$ of
the form $e \widetilde B e$, where $\widetilde B$ is a filtered
algebra such that $\gr \widetilde B \cong \cO_{\CC^{2n}} \rtimes
S_{n+1}$, and $e = \frac{1}{(n+1)!} \sum_{\sigma \in S_{n+1}} \sigma$
is the symmetrizer.
\begin{corollary} \label{defcor6}
  Let $B$ an arbitrary noncommutative spherical rational Cherednik
  algebra deforming $\cO_{\CC^{2n}}^{S_{n+1}}$.
  Then, $\dim \HH_0(B)= 1$.
\end{corollary}
In particular, this also gives another proof of the result from
\cite{BEG} that $B$ can have at most one irreducible
finite-dimensional representation.
\begin{remark}
  If $B$ admits any other filtered quantizations aside from the
  Cherednik algebras, then for these one concludes at least that $\dim
  \HH_0(B) \leq 1$ and $B$ admits at most one finite-dimensional
  irreducible representation. However, we do not know if there exist
  any quantizations other than the Cherednik algebras; cf.~the comments
  in \S \ref{ss:csr} below.
\end{remark}

\subsubsection{Prime ideals of quantizations}
Returning to the case of $S^n Y$ where $Y$ is affine symplectic, we
remark that there can never be any finite-dimensional representations
of quantizations of $S^n Y$ when $\dim Y > 0$ and $Y$ is connected,
since $S^n Y$ has no zero-dimensional symplectic leaves (i.e.,
subvarieties closed under the flow of Hamiltonian vector fields $\xi_f
:= \{f, -\}$). In more detail, recall that a primitive ideal of an
associative algebra $A$ is the kernel of an irreducible
representation.  If $A$ is a filtered quantization of an affine
Poisson variety $X$, then a primitive ideal $J$ is the kernel of a
finite-dimensional representation if and only if the support of $\gr
J$ is zero-dimensional.  
However, it is well known (and easy to check)
that the support of $\gr J$ must be closed under Hamiltonian flow on
$X$.  Since $S^n Y$ has no zero-dimensional symplectic leaves, it
follows that $A$ cannot have any finite-dimensional irreducible
representations.
% Alternatively, if one supposed that $\gr J$ were zero-dimensional,
% then one could complete $Y$ at the collection of points in $Z(\gr
% J)$, and obtain a nonzero representation of a quantization of the
% completed $\cO_{S^n \hat Y}$, contradicting the fact that $\hat Y$
% has no top cohomology.

However, we can still make a nontrivial statement about more general
primitive ideals of quantizations of $\cO_{S^n Y}$.  In fact, we can
consider more generally prime ideals: recall that a (two-sided) ideal
$J \subseteq A$ is prime if $R = A/J$ is a prime ring, i.e., $aRb = 0$
if and only if either $a$ or $b$ is zero.  All primitive ideals are
prime.

Using the method of I. Losev's appendix to \cite{ESdm}, we may then
deduce
\begin{corollary}\label{c:los1}
  Let $Y$ be connected affine symplectic 
  and let $A_\hbar$ or $B$ be a
  deformation or filtered quantization of $\cO_{S^n Y}$,
  respectively. For each $i \leq n$, the number of prime
  ideals of $A_\hbar[\hbar^{-1}]$ or $B$ (over $\CC\dblparens{\hbar}$ or
  $\CC$, respectively) whose support has codimension $i \dim Y$ in $S^n Y$
  is at most $p_{n,i}$, which is given by the generating function
\begin{equation}
\sum_{n, i \geq 0} p_{n,i} s^i t^n = \prod_{m \geq 0} \frac{1}{1-s^mt^{m+1}},
\end{equation}
i.e., the number of partitions of $n$ with $n-i$ parts.  
There are no prime ideals whose support has codimension not a multiple of $\dim Y$.
\end{corollary}
In particular, the bound on the number of prime ideals is independent of $Y$. 

Similarly, in the case of $V^n/S_{n+1}$, we may deduce
\begin{corollary}\label{c:los2}
  Let $V$ be a symplectic vector space, and $A_\hbar$ or $B$ be a
  deformation or filtered quantization of $\cO_{V^n}^{S_{n+1}}$. Then,
  for each $i \leq n$, the number of prime ideals of
  $A_\hbar[\hbar^{-1}]$ or $B$ whose support has codimension $i \dim
  V$ in $V^n / S_{n+1}$ is at most $p_{n+1,i}$.
\end{corollary}

\subsubsection{Poisson deformations and
zero-dimensional symplectic leaves}
Assume that $A$ is a graded algebra which is Poisson with Poisson
bracket of degree $-d$. Then, recall that a filtered Poisson
deformation is a Poisson algebra $B$ whose Poisson bracket satisfies
$\{B_{\leq i}, B_{\leq j}\} \subseteq B_{\leq i+j-d}$, and such that
$\gr B = A$ as a Poisson algebra.  Similarly, if $A$ is an arbitrary
Poisson algebra, one can consider Poisson algebras $(A_\hbar, \star,
\{-,-\}_\star)$ over $\CC\dblsqbrs{\hbar}$, which are isomorphic to
$A\dblsqbrs{\hbar}$ as $\CC\dblsqbrs{\hbar}$-modules, and satisfy $a \star b = ab +
O(\hbar)$ and $\{a, b\}_\star = \hbar \{a,b\} + O(\hbar^{2})$ for all
$a, b \in A$. Let us call these \emph{formal Poisson deformations}.
The two deformations above are analogous to filtered and deformation
quantizations, respectively (there is a slight discrepancy with the
use of the term ``deformation'' which refers to a formal parameter in
the quantization case but not in the Poisson case).

In the filtered case, one has a surjection $\HP_0(A) \onto \gr \HP_0(B)$,
and in the formal case, one has $\HP_0(A)\dblparens{\hbar} \onto \gr \HP_0(B_\hbar[\hbar^{-1}])$.  

Finally, recall that a zero-dimensional symplectic leaf of a Poisson
variety $X$ is a point $x \in X$ at which all Hamiltonian vector
fields vanish. Equivalently, $x$ is a point at which the evaluation
map $\ev_x: \cO_X \to \CC$ is a Poisson trace.  Note that the
evaluation maps at distinct points of $X$ are linearly independent.

Therefore, as before, we deduce
\begin{corollary}\label{c:fpoiss1}
Let $Y$ be a connected affine symplectic variety.
If $B_\hbar$ is a formal Poisson deformation of $S^n Y$, then $\dim_{\CC\dblparens{\hbar}} \HP_0(B_\hbar[\hbar^{-1}]) \leq a_n(\dim H^{\dim Y}(Y))$. 
\end{corollary}
Now consider the linear case.  By \cite[Proposition 1.16]{GiKal}, the
second Poisson cohomology $\HP^2(\cO_{V^n}^{S_{n+1}})$ is zero unless
$\dim V = 2$, in which case it is one-dimensional. So, there are only
nontrivial (filtered or deformation) Poisson deformations when $\dim V
= 2$. Moreover, by \cite[Theorem 1.18]{GiKal}, there is a universal
Poisson deformation of $\cO_{\CC^{2n}}^{S_{n+1}}$, given by the family
of commutative spherical rational Cherednik algebras. Using results of
\cite{EGsra}, we may conclude
\begin{corollary} \label{c:fpoiss2} Let $V$ be a symplectic vector
  space and $B$ be a nontrivial filtered Poisson deformation of
  $\cO_{V^n}^{S_{n+1}}$.  Then $V \cong \CC^2$, $B$ is a commutative
  spherical rational Cherednik algebra, $\dim \HP_0(B) = 1$, and
  $\Spec B$ has at most one zero-dimensional symplectic leaf.
\end{corollary}

\subsection{A canonical $\caD$-module on $S^n Y$ for $Y$ symplectic}
\label{ss:typea-int-dmod}
Here we will explain and generalize Theorem \ref{symphp0thm} using
$\caD$-modules.

We first recall the basic construction from \cite{ESdm} for Poisson
varieties.  Let $X$ be an affine Poisson variety, i.e., $X = \Spec A$
where $A$ is a Poisson algebra over $\CC$ which is finitely generated
as an algebra over $\CC$.  Let $i: X \into V$ be an embedding of $X$
into a smooth affine variety $V$.  Then, \cite{ESdm} defined the right
$\caD_V$-module $M(X,i)$ on $V$ as the quotient of the ring $\caD(V)$
of differential operators on $V$ with polynomial coefficients by the
right ideal generated by functions vanishing on $X$ and vector fields
which, on $X$, are parallel to $X$ and restrict to Hamiltonian vector
fields.  As explained there, this does not depend on the choice of
embedding $i: X \into V$, in the sense that, given two embeddings
$i_1: X \into V_1$ and $i_2: X \into V_2$, the resulting
$\caD_V$-modules $M(X,i_1)$ and $M(X,i_2)$ are images of each other
(up to isomorphism) under Kashiwara's equivalence of categories of
$\caD_V$-modules on $V_1$ and $V_2$ supported on $X$.  We may thus
refer to the module as $M(X)$ when not using the embedding.  (Note
that one can also define $M(X)$ without using an embedding at all, as
a quotient of the canonical right $\caD$-module $\caD(X)$ by the left
action of Hamiltonian vector fields: see \cite{ESdm}.)

% \footnote{One can define $M(X)$ without
%   using an embedding at all: although the correct category, call it
%   $\caD(X)$, of right $\caD$-modules on $X$ is not the same as the
%   category of modules over the ring $\caD_X$ or any other ring, it can
%   be intrinsically defined using, e.g., crystals, and it is equipped
%   with a canonical right $\caD$-module called $\caD_X$, such that
%   $\Hom_{\caD(X)}(\caD_X, N) = \Gamma(X,N)$, the global sections
%   functor. Also, $\caD_X$ is equipped with a canonical left action by
%   vector fields on $X$. Then, $M(X)$ is the quotient of $\caD_X$ by
%   the submodule generated by the left action of Hamiltonian vector
%   fields. See \cite{ESdm}.}

The motivation for the definition of $M(X)$ is the formula
$\pi_0(M(X)) \cong \HP_0(\cO_X)$, where $\pi: X \rightarrow \text{pt}$
is the projection to the point, and $\pi_0$ is the underived direct
image.

On the other hand, since the definition of the $\caD$-module $M(X)$ is local,
as explained in \cite{ESdm}, it makes sense to define $M(X)$ even in the
case that $X$ is not affine. 

We now present a theorem giving the structure of $M(X)$ when $X = S^n
Y$, for $Y$ a symplectic variety that need not be affine.  Let
$\Delta_i: Y \into S^i Y$ be the diagonal embedding, and for
$\sum_{j=1}^k r_j i_j = n$, let $q: (S^{i_1} Y)^{r_1} \times \cdots
\times (S^{i_k} Y)^{r_k} \onto S^n Y$ be the obvious projection.
% Let $q_*^{S_n}$ be the equivariant pushforward functor from
% $D$-modules on $Y^n$ to $S^n Y$, given by $(q_n)_*^{S_n}(M) :=
% (q_n)_*(M)^{S_n}$. This strengthens Theorems \ref{symphp0thm} and
% \ref{snp1thm}, as we will explain.
\begin{theorem}\label{sympdmthm}
 \begin{equation} \label{e:sympdmthm}
M(S^n Y) \cong \bigoplus_{r_1 \cdot i_1 + \cdots + r_k \cdot i_k = n, 1 \leq i_1 < \cdots < i_k, r_j \geq 1\, \forall j}
  q_* \bigl((\Delta_{i_1})_*(\Omega_Y)^{\boxtimes r_1} \boxtimes \cdots \boxtimes
  (\Delta_{i_k})_*(\Omega_Y)^{\boxtimes r_k} \bigr)^{S_{r_1} \times \cdots \times S_{r_k}}.
\end{equation}
\end{theorem}
\begin{remark} \label{r:Mphi} Let $Y$ be connected and $V$ be a
  symplectic vector space of dimension equal to the dimension of $Y$.
  The results of \cite{ReScmat} on $\HP_0(\cO_{V^{n-1}}^{S_{n}},
  \cO_{V^{n-1}})$ (cf.~Remark \ref{r:ReSc}) imply information on the
  structure of the $S_n$-equivariant $\caD$-module $M_{\phi}(Y^n)$,
  where $\phi: Y^n \onto Y^n/S_n = S^n Y$, and $M_\phi(Y^n)$ is
  defined as the quotient of $D_{Y^n}$ by the $S^n$-invariant
  Hamiltonian vector fields (more generally, for $\psi: X \to Z$,
  $M_\psi(X)$ is locally the quotient of $\caD_X$ by the action of
  Hamiltonian vector fields associated to Hamiltonians pulled back
  from $Z$: see \cite{ESdm}).  The computation of $M_\phi(Y^n)$
  reduces to a study of the diagonal $Y \into Y^n$ as before, and
  there one obtains the $S_n$-module $\Omega_Y \otimes
  \HP_0(\cO_{V^{n-1}}^{S_{n}}, \cO_{V^{n-1}})$.  Then, the results of
  \emph{op.~cit.} mentioned in Remark \ref{r:ReSc} say that the
  reflection representation $\mathfrak{h}$ of $S_n$ does not occur
  here, and the representation $\wedge^2 \mathfrak{h}$ occurs there
  tensored by $\lfloor \frac{n-1}{2} \rfloor$ copies of $\Omega_Y$.
  Using this and the theorem, one can obtain analogues of Theorem
  \ref{symphp0thm} giving information on the structure of
  $\HP_0(\cO_{S^n Y}, \cO_{Y^n})$.
\end{remark}
This implies the following ``derived'' generalization of Theorem
\ref{symphp0thm}.  For any affine Poisson variety $X$, let
$\HP^{DR}_i(X) := L^i \pi_*(M(X))$ be the $i$-th derived pushforward
of $M(X)$ to a point, where $\pi: X \rightarrow \pt$ is the
projection.  This is called the $i$-th \emph{Poisson-de Rham homology}
of $X$ and was defined in \cite{ESdm}. Note that $\HP^{DR}_0(X) =
\HP_0(\cO_X)$. Moreover, when $X$ is symplectic and connected,
$\HP^{DR}_i(X) \cong H^{\dim X - i}(X) \cong \HP_i(\cO_X)$, where the
latter is the usual $i$-th Poisson homology (which differs from
$\HP^{DR}_i$ for general affine $X$ when $i > 0$).
\begin{corollary}\label{symphp0cor}
  Let $Y$ be affine symplectic and connected.
  Then, as bigraded algebras (in de Rham degree and in symmetric power plus
  $t$ degrees, i.e., $|\HP^{DR}_i(S^n Y)| = (i,n)$ and $|t|=(0,1)$),  
\begin{equation}\label{hpdreqn1}
  \bigoplus_{n \geq 0} \HP^{DR}_\bullet(S^n Y)^* \cong \Sym(\HP^{DR}_\bullet(Y)^*[t]) = \Sym(H^{\dim Y - \bullet}(Y)^*[t]).
\end{equation}
\end{corollary}
Next, continuing to assume that $Y$ is affine symplectic and connected,
let $A_\hbar$ be a deformation quantization of $\cO_{Y}$. Then, we deduce the following generalization of Corollary \ref{defcor1}:
\begin{corollary}\label{c:defcor1-dr} Taking the $\C\dblparens{\hbar}$-linear dual, we
obtain an isomorphism of bigraded algebras over $\C\dblparens{\hbar}$ (with $|\HH_i(\Sym^n A_\hbar)[\hbar^{-1}]^*| = (i,n) = |\HP_i^{DR}(S^n Y)\dblparens{\hbar}^*|$):
\begin{equation} \label{e:defcor1-dr}
\bigoplus_{n \geq 0} \HH_\bullet(\Sym^n A_\hbar)[\hbar^{-1}]^* \cong
\bigoplus_{n \geq 0} \HP_\bullet^{DR}(S^n Y)\dblparens{\hbar}^*.
\end{equation}
\end{corollary}

\subsection{Conjectures on symplectic resolutions}\label{ss:csr}
In this subsection, we explain some conjectures related to symplectic
resolutions motivated by the preceding results and also
\cite{ESsym}. The material of this subsection will not be needed
elsewhere in this paper.
 
For an irreducible (affine) Poisson variety $X$, we say that a morphism
$\widetilde{X} \to X$ is a symplectic resolution if $\widetilde{X}$ is
symplectic and the morphism is proper, birational, and Poisson (the
latter condition means that its pullback is a morphism of sheaves of
Poisson algebras). 
% Note that, in this case, $\cO_X = \Gamma(\widetilde{X},
% \cO_{\widetilde{X}})$ as Poisson algebras.

When $Y$ is a connected affine symplectic surface, $S^n Y$ admits a
symplectic resolution $\Hilb^n Y \onto S^n Y$, and we can deduce
from Corollary \ref{symphp0cor} and the known description of the cohomology of $\Hilb^n Y$ that
$\HP^{DR}_\bullet(S^n Y) \cong H^{2n - \bullet}(\Hilb^n Y)$.
  This suggests
\begin{conjecture}\label{c:hpresconj}
  Let $X$ be an irreducible affine Poisson variety with 
  a symplectic resolution $\rho: \widetilde{X} \onto X$. Then:
\begin{enumerate}
\item[(a)] $\HP_0(\cO_X) \cong H^{\dim \widetilde{X}}(\widetilde X)$.
\item[(b)] $\HP^{DR}_\bullet(X) \cong H^{\dim \widetilde{X} - \bullet}(\widetilde X)$.  
\item[(c)] $M(X) \cong \rho_* \Omega_{\widetilde X}$.
\end{enumerate}
\end{conjecture}
In part (c), $\rho_*$ refers to the derived pushforward.  Clearly,
$(c) \Rightarrow (b) \Rightarrow (a)$. Also, we remark that part (c)
makes sense even when $X$ is not affine, so it is reasonable to
conjecture that the affine assumption is not needed (and this also
would imply the generalization of (b) to nonaffine $X$, if we extend the
definition of $\HP^{DR}_\bullet(X)$ by taking the appropriate derived
pushforward of $M(X)$ to a point).

Note that (c) would imply that $M(X)$ is semisimple holonomic with
regular singularities, by the decomposition theorem \cite[Th\'eor\`eme
6.2.5]{BBDfp} (although the holonomicity already follows from
\cite[Theorem 3.1]{ESdm} once one notices that $X$ necessarily has
finitely many symplectic leaves; however, as pointed out in \cite[Example
4.11]{ESdm}, the latter condition does not imply that the
singularities are regular, and in fact neither does it imply that $M(X)$ is
semisimple). Similarly, it would follow immediately
from the conjecture that $\rho_* \Omega_{\widetilde X}$ is a
$\caD$-module rather than a complex, although this already follows
from the fact, \cite[Lemma 2.11]{Kalss}, that $\rho$ is a semismall
morphism.

We can prove the conjecture in three
cases:
\begin{enumerate}
\item[(A)] If $\widetilde X = \Hilb^n Y$ and $X = S^n Y$, part (c)
  follows from Theorem \ref{sympdmthm} together with the standard
  computation of $\rho_* \Omega_{\Hilb^n Y}$ (see \cite[Theorem
  3]{GoSohilb}).
\item[(B)] If $\widetilde X = T^* (G/B)$ is the Springer resolution of
  the nilpotent cone $X \subseteq \operatorname{Lie} G$, for $G$ a
  semisimple connected complex Lie group and $B < G$ a Borel, or more
  generally the restriction of this to the resolution of a Slodowy
  slice of $X$ (a transverse slice at a point $e \in X$ to its
  coadjoint orbit), part (c) follows from the main result of
  \cite{ESwalg}.
% (which is based on a computation of Hotta and Kashiwara).
\item[(C)] If $\displaystyle \widetilde X = \Hilb^n (\widetilde
  {\CC^2/G}) \mathop{\onto}^{\rho} \Sym^n(\CC^2/G) = X$, where $G <
  \SL_2(\CC)$ is a finite subgroup, i.e., $X$ is a symmetric power of
  a Kleinian singularity, and $\widetilde {\CC^2/G}$ is the minimal
  resolution of the Kleinian singularity. Then, the argument is
  similar to that of (i), using the computation of $\HP_0(\cO_X)$ from
  \cite[Theorem 1.1.14]{ESsym}: see \S \ref{s:conj-proof-klein} below.
\end{enumerate}

\begin{remark}\label{r:conj-gr}
  We stress that, for all parts (a)--(c) of the conjecture, 
  we only conjecture an abstract
  isomorphism (which is confirmed in the cases mentioned above),
  not a canonical isomorphism; i.e., in the cases of (a) and (b), we
  conjecture only an equality of dimensions. It would be desirable to
  refine the conjecture to give a more precise relationship between
  the two conjecturally isomorphic objects.

  At least for part (a), we can do this: 
  \eqref{e:hh-quant-fiber} below should imply that, for
  suitable deformation quantizations $B_\hbar$ of $X$, one has a \emph{canonical}
  isomorphism $\HH_{\bullet}(B_\hbar[\hbar^{-1}]) \cong H^{\dim X -
    \bullet}(\widetilde X)\dblparens{\hbar}$.  Since there is a canonical
  surjection $\HP_0(\cO_X)\dblparens{\hbar} \onto \gr \HH_0(B_\hbar[\hbar^{-1}])$,
  % spectral sequence $\HP_{\bullet}(\cO_X)\dblparens{\hbar}
  % \Rightarrow \HH_{\bullet}(B_\hbar[\hbar^{-1}])$,
  this suggests that, in the formal version with $\hbar$, there may
  rather be a filtration on the right-hand side of (a) whose
  associated graded vector space is the left-hand side, i.e., that
  there is a canonical isomorphism $\HP_0(\cO_X)\dblparens{\hbar} \iso
  \gr H^{\dim X}(\widetilde X)\dblparens{\hbar}$.
    
  Moreover, if $X$ has a contracting $\CC^\times$ action, we can
  eliminate the $\hbar$ using this grading, and should obtain a
  canonical isomorphism $\HP_0(\cO_X) \iso \gr H^{\dim X}(\widetilde
  X)$.  This holds in all cases we have checked (e.g., cases (B) and
  (C); note that there is no $\CC^\times$ action in case (A) in
  general).

  For parts (b) and (c) of the conjecture, it would be desirable to
  have a similar statement. However, we know of no direct relationship
  between the Poisson-de Rham homology of $X$ and the Hochschild
  homology of a quantization (there is a spectral sequence from
  ordinary Poisson homology of $\cO_X$ to this Hochschild homology,
  but ordinary Poisson homology only coincides with Poisson-de Rham
  homology in degree zero).  Perhaps this problem could be alleviated
  using the universal formal deformation $\mathcal{X}$ of $\widetilde
  X$ of \cite[Theorem 1.1]{KVpmnchsm} discussed below, which is
  generically affine symplectic, and which maps to the formal
  deformation $\Spec \Gamma(\mathcal{X},\cO_{\mathcal{X}})$ of $X$ in
  a way which is generically an isomorphism, since for affine
  symplectic varieties the Poisson-de Rham and ordinary Poisson
  homology coincide.
\end{remark}

 Next, we can pose a conjecture on the Hochschild homology of
 quantizations.  To motivate this, note that, in the case of (A) above, if 
  $A_\hbar$ is a deformation quantization of $\cO_Y$, then 
  % both $\HP_\bullet^{DR}(X)$ and $H^{\dim X-\bullet}(\widetilde X)$
  % have the same dimension as $\HH_\bullet(\Sym^n
  % A_\hbar[\hbar^{-1}])$. Hence,
by Corollary \ref{c:defcor1-dr},
$\HP^{DR}_\bullet(S^n Y)\dblparens{\hbar} \cong 
\HH_\bullet(\Sym^n A_\hbar[\hbar^{-1}])$.  We would like to generalize
this to the case of general symplectic resolutions.

 We will be particularly interested in quantizations obtainable by
 quantizing the symplectic resolution in the sense of \cite{BKfqac}.
 Namely, according to \cite[Definition 1.3]{BKfqac}, a quantization of
 $\widetilde X$ is a sheaf $\mathcal{B}_\hbar$ of associative flat
 $\CC\dblsqbrs{\hbar}$-algebras on $X$ equipped with an isomorphism
 $\mathcal{B}_\hbar/\hbar \mathcal{B}_\hbar \cong \cO_{\widetilde X}$.
 We will additionally require that the induced Poisson structure on
 $\cO_{\widetilde X}$ is the one coming from the symplectic form.  By
 \cite[Theorem 1.8]{BKfqac}, there is a semiuniversal family of such
 quantizations, parameterized by $\hbar H^2(\widetilde X)\dblsqbrs{\hbar}$.
 (Moreover, it seems reasonable to ask if these produce all
 quantizations of $X$, or if there is a semiuniversal family of all
 quantizations in which these map to a dense subset.)
\begin{conjecture}\label{c:hpquantconj}
  Let $X$ be an irreducible affine Poisson variety which admits a symplectic
  resolution.
\begin{enumerate}
\item[(i)] For every deformation quantization $A_\hbar$ of $\cO_X$,
  % obtained as the global sections of a quantization of $\widetilde
  % X$,
  the canonical surjection is an isomorphism
\begin{equation}\label{e:hpquantconj}
  \HP_0(\cO_X)\dblparens{\hbar} \iso \gr
  \HH_0(A_\hbar[\hbar^{-1}]).
\end{equation}
\item[(ii)] There is a countable collection of $\hbar$-homogeneous
hypersurfaces in $\hbar H^2(\widetilde X)\dblsqbrs{\hbar}$ 
such that, for $A_\hbar$ obtained as the global sections of
a quantization in the family $\hbar H^2(\widetilde X)\dblsqbrs{\hbar}$ outside of this collection, one has (abstractly)
\begin{equation} \label{e:hpdrquantconj}
\HP_\bullet^{DR}(X)\dblparens{\hbar} \cong \gr \HH_\bullet(A_\hbar[\hbar^{-1}]).
\end{equation}
% \item[(iii)] For every deformation quantization $\mathcal{B}_\hbar$
%   of $\widetilde X$,
%   \begin{equation}
% \HP_*^{DR}(X)\dblparens{\hbar} \cong \gr \HH_*(\mathcal{B}^{\CC\dblparens{\hbar}}_\hbar[\hbar^{-1}]).
% \end{equation}
\end{enumerate}
\end{conjecture}
Here, by an $\hbar$-homogeneous hypersurface in $\hbar 
H^2(\widetilde X)\dblsqbrs{\hbar}$, we mean by definition a subvariety
of the form $Z \times \hbar^{m}H^2(\widetilde X)\dblsqbrs{\hbar}$,
where $Z \subseteq \bigoplus_{j=1}^{m-1} \hbar^j H^2(\widetilde
X)$ is cut out by an equation which is homogeneous in $\hbar$ of some
degree.
\begin{remark}
  As in Remark \ref{r:conj-gr} above, it would be better if in (ii)
  one could construct a canonical map from the LHS to the RHS which is
  conjecturally an isomorphism, but we are not sure how to do this.
\end{remark}
Moreover, given a semiuniversal quantization of $X$, one can ask if
\eqref{e:hpdrquantconj} still holds for this family. Note that (ii)
implies (i) (for quantizations considered in (ii)), since $\dim_{\CC\dblparens{\hbar}}
\HH_0(A_\hbar[\hbar^{-1}])$ is upper semicontinuous and bounded above
by $\dim \HP_0(\cO_X)$ (and $\HP_0(\cO_X) = \HP_0^{DR}(X)$, unlike in
higher degrees).

Also, note that the genericity assumption of (ii) above is needed:
already in the case $X = \CC^2/(\ZZ/2)$, there exist quantizations for
which \eqref{e:hpdrquantconj} does not hold (this follows from
\cite[Theorem 2.1]{FSSAhhc}; see also \cite[Remark 1.14]{ESwalg}).
Indeed, only in degree zero does one obtain a (natural) surjection
from $\HP_\bullet^{DR}(X)\dblparens{\hbar}$ to
$\gr \HH_\bullet(A_\hbar[\hbar^{-1}])$.

In cases (B) and (C) above, we can prove this conjecture, at least
when (i) is restricted to quantizations coming from the symplectic
resolution.  In case (B), one should be able to check that the
algebras $A_\hbar$ appearing in the conjecture are the Rees algebras
of the quantum $W$-algebras deforming $\cO_X$.
% (one can even explicitly obtain these using the explicit
% quantization of the Grothendieck-Springer resolution constructed by
% Dodd-Kremnizer and Ginzburg).
For these algebras, parts (i) and (ii) of the conjecture follow from
\cite[Theorem 1.10.(ii)]{ESwalg} and \cite[Theorem 1.13]{ESwalg}.

In case (C), the algebras $A_\hbar$ appearing in the conjecture should
be the Rees algebras of the spherical symplectic reflection algebras
\cite{EGsra} deforming $\cO_{\Sym^n (\CC^2/G)} = \cO_{\CC^{2n}}^{G^n
  \rtimes S_n}$.
% (in the case $G$ is trivial, one can explicitly use the quantization
% of $\Hilb^n \CC^2$ constructed by Gordon and Stafford)
For these algebras, part (i) of the conjecture
is a consequence of \cite[Corollary 1.3.2]{ESsym}. Part (ii) follows
by comparing the explicit description of $M(X)$ given in \S
\ref{s:conj-proof-klein} below (for the LHS) with the description of
$\HH_\bullet(\caD_{\CC^n}^{G^n \rtimes S_n})$ from \cite{AFLS}, as
well as the fact from \cite[Theorem 1.8]{EGsra} that this coincides
with $\HH_\bullet(A)$ for generic spherical symplectic reflection
algebras $A$ quantizing $\cO_{\CC^{2n}}^{G^n \rtimes S_n} =
\cO_{\Sym^n (\CC^2/G)}$.

If true, the conjecture would yield a necessary criterion for
existence of symplectic resolutions (where in (ii) we take a
semiuniversal family of quantizations). This condition does \emph{not}
appear to be sufficient, however: already in the case that $X = \Sym^n
V$ for $V$ a symplectic vector space of dimension $\geq 4$, our main
theorem implies that \eqref{e:hpdrquantconj} holds for the
quantization $\Sym^n \Weyl(V)$. We are not sure if there exist other
quantizations: for the quasiclassical analogue, there exist no
nontrivial Poisson deformations as discussed after Corollary
\ref{c:fpoiss1}. On the other hand, $X$ does not admit a symplectic
resolution by \cite{Verhsgos} (since $G$ is not generated by
symplectic reflections, i.e., elements $g \in G$ such that $g-\Id$ has
rank two; in fact, $G$ has no symplectic reflections, which is why
$\HP^2(\cO_X) = 0$).

Finally, we remark that Conjecture \ref{c:hpresconj} almost implies
Conjecture \ref{c:hpquantconj} (at least if we restrict part (i) to
quantizations coming from the resolution).  First of all, by
\cite[Theorem 1.1]{KVpmnchsm}, there is a universal formal deformation
$\mathcal{X}$ of $\widetilde X$ in the category of symplectic schemes,
which lies over the formal completion $\widehat H^2(\widetilde X)$ of
$H^2(\widetilde X)$ at the origin. By \cite[Theorem 1.8, Lemma
6.4]{BKfqac}, $\mathcal{X}$ also admits a canonical quantization over
$\widehat H^2(\widetilde X)$, so that the quantization
$\mathcal{B}_\hbar$ corresponding to a formal power series $P \in
\hbar H^2(\widetilde X)\dblsqbrs{\hbar}$ is the pullback of the canonical
quantization $\mathcal{B}_{\mathcal{X}}$ of $\mathcal{X}$ by the formal point
$p \in \widehat H^2(\widetilde X)$ corresponding to $P$. Now,
according to \cite[Lemma 2.5]{Kal-deq}, for generic $p$,
the fiber of $\mathcal{X}$
over $p$ is affine. For such $p$, it should follow that
\begin{equation} \label{e:hh-quant-fiber}
\HH_*(\Gamma(\mathcal{X}_p, p^*
\mathcal{B}_{\mathcal{X}}[\hbar^{-1}])) \cong H^{\dim X -
  *}(\mathcal{X}_p)\dblparens{\hbar} \cong
H^{\dim X - *}(\widetilde
X)\dblparens{\hbar},
\end{equation}
% using affineness for the first isomorphism and Conjecture
% \ref{c:hpresconj}.(b), and
adapting the usual identification of Hochschild homology of
quantizations of an affine symplectic variety with the de Rham
cohomology of the variety for the first isomorphism, and applying
topological triviality of the family of deformations for the second
isomorphism.  This would yield Conjecture \ref{c:hpquantconj}.(ii).
Then, to deduce part (i), we apply part (ii) together with the fact
that $\dim \HP_0(\cO_X) \geq \dim_{\CC\dblparens{\hbar}} \HH_0(A_\hbar[\hbar^{-1}])$ for all
quantizations, along with upper semicontinuity of $\dim
\HH_0(\Gamma(\mathcal{X}_p, p^*
\mathcal{B}_{\mathcal{X}}[\hbar^{-1}]))$ in $p$.  Thus, Conjecture
\ref{c:hpresconj} should also imply Conjecture \ref{c:hpquantconj}, at
least in (ii) if we ask only for an abstract isomorphism of
$\CC\dblparens{\hbar}$-vector spaces which preserves the homological grading
($\bullet$).

\subsubsection{The case of linear quotient singularities}
In the case when $X = V/G$ is a linear quotient singularity with $G <
\Sp(V)$, the main result of \cite{AFLS} computes $\dim
\HH_{2i}(\Weyl(V)^G)$: this is the number of conjugacy classes of $g
\in G$ such that $\dim \ker(g-\Id) = 2i$. Here, $\Weyl(V)$ is the Weyl
algebra and $\Weyl(V)^G$ is therefore a filtered quantization of
$\cO_V^G$.  This would imply the first part of the following
conjecture:
\begin{conjecture}\label{c:hpdim}
  Suppose that $G < \Sp(V)$ is finite and $V/G$ admits a symplectic
  resolution.  Then 
\begin{enumerate}
\item[(i)] The canonical surjection is an isomorphism
$\HP_{0}(\cO_V^G) \iso \gr \HH_0(\Weyl(V)^G)$.  In particular, $\dim \HP_0(\cO_V^G)$ is
the number of conjugacy classes of
  elements $g \in G$ such that $g - \Id$ is invertible.
\item[(ii)] For all $i \geq 0$, abstractly, $\HP_{2i}^{DR}(V/G) \cong
  \HH_{2i}(\Weyl(V)^G)$, i.e., $\dim \HP_{2i}^{DR}(V/G)$ is the number
  of conjugacy classes of $g \in G$ such that $\ker(g-\Id)$
  has dimension $2i$ (and $\HP_{2i+1}^{DR}(V/G) = 0$).
\end{enumerate}
\end{conjecture}
Conversely, part (i) of the above conjecture would imply
Conjecture \ref{c:hpquantconj} for noncommutative spherical symplectic
reflection algebras deforming
$\cO_V^G$: this follows from reasoning similar to the proof of
Corollary \ref{defcor6}.  Namely, by \cite[Theorem 1.8]{EGsra}, for
generic such algebras, $\dim \HH_0(A)$ coincides with $\dim
\HH_0(\Weyl(V)^G)$.  Hence, by upper semicontinuity of $\dim \HH_0$ in
the family, $\dim \HH_0(A)$ is at least $\dim \HH_0(\Weyl(V)^G)$ for
all noncommutative spherical symplectic reflection algebras 
deforming $\cO_V^G$.  Then, Conjecture
\ref{c:hpdim} would say that this coincides also with $\dim
\HP_0(\cO_V^G)$, which is an upper bound for $\dim \HH_0(A)$ for all
quantizations $A$. Hence, $\dim \HH_0(A)$ would be constant in the
family, implying Conjecture \ref{c:hpquantconj}.(i) for such algebras.

In case (C) above, i.e., for $G$ a wreath product of a finite subgroup
of $\SL_2(\CC)$ with $S_n$ for some $n \geq 1$, we can prove the above
conjecture: it follows from our proof of Conjecture \ref{c:hpresconj}
below (or alternatively, it follows from Conjecture
\ref{c:hpquantconj}, since in this case the family of quantizations
obtained from the resolution of singularities is exactly the
noncommutative spherical symplectic reflection algebras).  Note that,
in this case, statement (i) was a conjecture by J. Alev of
\cite[Remark 40]{But} (he possibly conjectured it for some
other groups $G$ as well elsewhere); this conjecture was
first proved in \cite{ESsym}, apart from the cases $n = 2$ and $n=3$
where it was proved in \cite{AlFo} and \cite{But}, respectively.

\subsubsection{Proof of Conjecture \ref{c:hpresconj} in the case $X =
  \Sym^n (\CC^2/G)$}\label{s:conj-proof-klein}
Let $X = \Sym^n (\CC^2/G)$ where $G < \SL_2(\CC)$ is a finite group.  By
\cite[Corollary 4.16]{ESdm}, $M(X)$ is a direct sum of IC $D$-modules
of the symplectic leaves with some multiplicities. These leaves are indexed
by tuples $(r, r_1, \ldots, r_k)$ of nonnegative integers, such that
$r+\sum_{j=1}^k j \cdot r_j = n$.  This symplectic leaf, $X_{(r,r_1,
  \ldots, r_k)}$, has closure 
given by the 
 image of
\begin{multline}
\{0\} \times \Sym^{r_1}(\CC^2/G) \times \cdots \times \Sym^{r_k}(\CC^2/G) \\ 
\into \{0\}^r \times \Sym^{r_1}((\CC^2/G)^{1}) \times \Sym^{r_2}((\CC^2/G)^2) \times \cdots \times \Sym^{r_k} ((\CC^2/G)^{k}) \to \Sym^n (\CC^2/G) = X.
\end{multline}
%(This is similar to \S \ref{ss:typed-int-dmod}). 
 The multiplicity of $IC(\overline{X_{(r,r_1,\ldots,r_k)}})$, by \emph{op.~cit.}, is $\dim \HP_0(\cO_{Z_{(r,r_1,\ldots,r_k)}})$, where
$Z_{(r,r_1,\ldots,r_k)} = \Sym^r (\CC^2/G) \times \prod_{j=2}^k (\CC^{2(j-1)}/S_j)^{r_j}$.
By Theorem \ref{snp1thm} and \cite[Theorem 1.1.14]{ESsym}, this multiplicity
is equal to the number of $\dim \HP_0(\CC^2/G)$-multipartitions of $r$.  By, e.g., \cite{AL}, $\dim \HP_0(\CC^2/G)$ is the number of isomorphism classes of nontrivial representations of $G$, which is well known to be the number of irreducible components of the fiber $\pi^{-1}(0)$ of the resolution of Kleinian singularities, $\pi: \widetilde {\CC^2/G} \to \CC^2/G$.  

Then, it remains to show that the above is the same as
$\rho_*\Omega_{\widetilde X}$. 
%  We can decompose $\rho$ as
% \[
% \widetilde X = \Hilb^n \widetilde{\CC^2/G} \mathop{\to}^{\rho_1}
% \Sym^n \widetilde (\CC^2/G) \mathop{\to}^{\rho_2} \Sym^n (\CC^2/G) =
% X.
% \]
We can argue similarly to the aforementioned result, \cite[Theorem
3]{GoSohilb} (which dealt with the case where the fibers of $\rho$
were irreducible). Namely, since the above map is semismall (this is
well known in this case, and is also more generally true for all
symplectic resolutions by the aforementioned \cite[Lemma
2.11]{Kalss}), $\rho_*\Omega_{\widetilde X}$ decomposes as a direct
sum of intermediate extensions of local systems (i.e., $\cO$-coherent
$\caD$-modules) on each symplectic leaf of $X$. Moreover, the local
systems occurring on each symplectic leaf are the top cohomology of
the fibers of $\rho$ restricted to that leaf. By restricting to a
formal neighborhood of a symplectic leaf, using the explicit
description of the symplectic leaves above, the computation reduces to
the case of the point $\{0\} \in \Sym^{n'} (\CC^2/G)$ for all $n' \leq
n$. In this case, we evidently get a direct sum of delta-function
$\caD$-modules, with multiplicity given by the number of irreducible
components of $\rho^{-1}(0)$ of dimension $n'$.  This is equal to the
number of $m$-multipartitions of $n'$, where $m$ is the number of
irreducible components of the zero fiber of $\widetilde{\CC^2/G} \to
\CC^2/G$.  This is, however, the same multiplicity as for $M(X)$, as
mentioned above. We conclude that $\rho_* \Omega_{\widetilde X} \cong
M(X)$, as desired.

\subsection{Examples of nontrivial local systems in $M(X)$}
Note that, in all of the examples of affine Poisson varieties $X$
studied thus far in this paper, $M(X)$ is a direct sum of intermediate
extensions of trivial local systems on the symplectic leaves of
$X$. \emph{Here and below, ``local system'' refers to an
  $\cO$-coherent $\caD$-module on a smooth variety.} Note that these
were all examples of the form $X = U/G$ with $U$ affine symplectic and
$G$ a finite group of symplectic automorphisms of $U$.  In this
subsection, which will not be required in the remainder of the paper,
we construct other examples of this form such that nontrivial local
systems do appear in $M(X)$. This fulfills the promise of
\cite[footnote 6]{ESdm}.

In fact, by \cite[Corollary 4.16]{ESdm}, $M(X)$ is always semisimple
if $X = V / G$ for $V$ a symplectic vector space and $G < \Sp(V)$
finite. Also, by \cite[Theorem 4.21]{ESdm}, whenever $X = U/G$, $U$ is
a symplectic variety (not necessarily a vector space or even affine),
and $G$ is a finite group of symplectic automorphisms of $U$, then
$M(X)$ is always a direct sum of intermediate extensions of
one-dimensional local systems on symplectic leaves of $X$; these local
systems all have monodromy valued in $\pm 1$. Moreover, there is a
simple necessary (but not sufficient) criterion for the local systems
to be nontrivial: roughly, the action of $G$ on normal bundles to
preimages of symplectic leaves must contain quaternionic
representations.  More precisely, let $X_0 \subseteq X$ be a
symplectic leaf, and fix $x \in X_0$ with preimage $u \in U$. Then,
there can only be a nontrivial local system appearing in $M(X)|_{X_0}$
if the $\Stab_{G}(u)$-representation $(T_u U)^\perp$ contains a
quaternionic irreducible representation.  In particular, this implies
the aforementioned result (which we also explain directly in \S
\ref{ss:sympdmthm-pf} below) that $M(S^n Y)$ is a direct sum of
intermediate extensions of trivial local systems, when $Y$ is a
symplectic variety.  This is because the $\Stab_G(u)$-representations
$(T_u U)^\perp$ are all products of (reducible) representations
$\CC^{2m}$ of $S_{m+1}$ associated to type $A_m$ Weyl groups, and in
particular all irreducible summands are of real, not quaternionic,
type.

Now, let $X = U/G$ where $U$ is a symplectic variety and $G$ is a
finite group of symplectic automorphisms. Let $x \in X_0$ and $u \in
U$ be as above.  Let us describe the local system $M(X)|_{X_0}$ more
explicitly.  As observed in \cite[\S 4]{ESdm}, this local system has
fiber $\HP_0(\cO_{(T_u U)^\perp}^{\Stab_G(u)})$.  The monodromy is
given by the composition $\pi_1(X_0) \to \Sp_{\Stab_G(U)}((T_u
U)^\perp) \to \Aut(\HP_0(\cO_{(T_u U)^\perp}))$, where
$\Sp_{\Stab_G(u)}((T_u U)^\perp)$ denotes the group of automorphisms
of the symplectic vector space $(T_u U)^\perp$ preserving the
$G$-action.  The first map is given by the Hamiltonian flow along
$X_0$, as explained in \emph{op.~cit.}.  Moreover, as explained in
\emph{op.~cit.}, since $M(X)$ is locally constant along Hamiltonian
vector fields, the first map factors through $\Sp_{\Stab_G(u)}((T_u
U)^\perp)/\Sp_{\Stab_G(u)}((T_u U)^\perp)^\circ \cong \prod_{Q \in
  R_q((T_u U)^\perp)} \ZZ/2$, where $R_q((T_u U)^\perp)$ denotes the
set of isomorphism classes of quaternionic representations of $\Stab_G(u)$
occurring in $(T_u U)^\perp$.

We therefore have to consider the two resulting maps: (a) $\pi_1(X_0)
\to \prod_{Q \in R_q((T_u U)^\perp)} \ZZ/2$, and (b) $\prod_{Q \in
  R_q((T_u U)^\perp)} \ZZ/2 \to \Aut(\HP_0(\cO_{(T_u U)^\perp}))$.  We
first consider (a):
\begin{claim} \label{cl:monodromy} For any symplectic vector space
  $V$, finite subgroup $G < \Sp(V)$ such that $V^G = \{0\}$,
  symplectic variety $Y$, and homomorphism $\pi_1(Y) \to \prod_{Q \in
    R_q(U)} \ZZ/2$, one can construct a symplectic variety $U$ with an
  action of $G$ such that
\begin{enumerate}
\item[(i)] $U^G \cong Y$;
\item[(ii)] For $u \in U^G$, $(T_u U)^\perp \cong V$ as symplectic 
$G$-representations;
\item[(iii)] The map  $\pi_1(U^G)
\to \prod_{Q \in R_q((T_u U)^\perp)} \ZZ/2$ coincides with the given map under (i) and (ii).
\end{enumerate}
\end{claim}
Using the claim, it will remain only to exhibit a pair $(V,G)$ such
that the map (b) is nonzero.  Let us explain such an example.  We
begin by describing the map (b) more explicitly.  It is easy to see
that each generator $1_Q \in \ZZ/2$ corresponding to $Q \in R_q((T_u
U)^\perp)$ maps to $(-\Id)^{\mu_Q}$, where $\mu_Q$ is the operator $f
\mapsto |f|_Q$, assigning to functions their parity of degree in any
fixed summand of $(T_u U)^\perp$ isomorphic to $Q$ (this parity of
degree is independent of the choice of summand).  In particular, if
$(T_u U)^{\perp}$ is itself an irreducible quaternionic
representation, $\mu_Q$ is the parity of the polynomial degree.

More generally, if $(T_u U)^\perp$ is a direct sum of distinct
irreducible quaternionic representations, then the image of $(1,
\ldots, 1) \in \prod_{Q \in R_q((T_u U)^\perp)} \ZZ/2$ in
$\Aut(\HP_0(\cO_{(T_u U)^\perp}))$ is $(-\Id)^{\deg}$, where $\deg$ is
the polynomial degree. In particular, this is nontrivial in the case
that $\HP_0(\cO_{(T_u U)^\perp}^{\Stab_G(u)})$ is nontrivial in odd
degrees.  

An example of a pair $(V, G)$ of a symplectic vector space $V$ and a
finite subgroup $G < \Sp(V)$ such that $\HP_0(\cO_{V}^G)$ is
nontrivial in odd degrees was exhibited in the appendix to
\cite{hp0bounds}: there $V = V_1 \oplus V_2 \oplus V_3$ with $V_i$
irreducible quaternionic representations of $G$, with $\dim V_1 = \dim
V_2 = \dim V_3 = 2^m$ for $m \geq 2$.  In particular, the smallest
dimension of such $V$ found there is $12$. We can apply this to the
claim with $Y = \CC^\times \times \CC$, where $\CC^\times$ is the
punctured complex plane, together with the map $\pi_1(Y) \cong \ZZ \ni
1 \mapsto (1, \ldots, 1)$.  The resulting $X$ has dimension $14$
(for the case $m=2$), and $M(X)|_{X_0}$ is nontrivial.
\begin{remark}
  The same analysis as above can be applied more generally to the
  $\caD$-module $M_\phi(U)$, where $\phi: U \to U/G$ is the quotient
  map (cf.~Remark \ref{r:Mphi}).  The only difference is that the
  fiber $\HP_0(\cO_{(T_u U)^\perp}^{\Stab_G(u)})$ is replaced by the
  $\Stab_G(u)$-representation $\HP_0(\cO_{(T_u U)^\perp}^{\Stab_G(u)},
  \cO_{(T_u U)^\perp})$.  Then, one produces examples of nontrivial
  local systems in $M_\phi(U)$ from any triple $(V, G, Y)$ as in the
  claim such that $H_1(Y, \ZZ/2) \neq 0$ with $V$ an
  \emph{irreducible} quaternionic representation of $G$, since in this
  case, $\HP_0(\cO_V^G, \cO_V)$ is already nontrivial in degree one,
  where it is $V$ itself.  (Here, by nontrivial, we mean that they are
  nontrivial even considered as ordinary local systems, not merely as
  $G$-equivariant local systems.)  For example, one can take $Y =
  \CC^\times \times \CC$, $V = \CC^2$, and $G < \SL_2(\CC)$ any finite
  nonabelian subgroup. Then, $M_\phi(U)$ is nontrivial, and $\dim U =
  4$. (Note that this is the minimum possible dimension of a
  symplectic variety $U$ such that, for some finite group of
  automorphisms $G$, $M_\phi(U)$ can restrict to a nontrivial local
  system on some some locally closed subvariety, which we may assume
  is the locus $\{u \in U: \Stab(u)=K\}$ for some subgroup $K < G$.)

  Finally, we can generalize the above argument to obtain information
  about the $G$-isotypic components of $M_\phi(U)$.  In particular, if
  $V$ is a direct sum of distinct irreducible quaternionic
  representations of $G$, and one constructs the associated
  $G$-variety $U$ as above, then any irreducible representation of $G$
  that occurs in odd degree in $\HP_0(\cO_V^G, \cO_V)$ also occurs
  tensored by a nontrivial local system in $M_\phi(U)|_{U^G}$.  

  For example, in the case $V = \CC^2$ and $G < \SL_2(\CC)$ is
  nonabelian, such irreducible representations of $G$ are exactly the
  ones occurring in the odd tensor powers of $V$, since these are the ones
  where $-\Id \in G$ acts by multiplication by $-1$. Moreover, for such
  irreducible representations of $G$, the isotypic part of
  $M_\phi(U)|_{U^G}$ occurs without a summand of the trivial local
  system on $U^G$.
\end{remark}
\begin{proof}[Proof of Claim \ref{cl:monodromy}]
We recall first the description given in
\emph{op.~cit.}: for each $Q \in R_q((T_u U)^\perp)$, let $((T_u
U)^\perp_Q \subseteq (T_u U)^\perp$ be the isotypic component of $Q$.
Note that 
\[\Sp_{\Stab_G(u)}((T_u U)^\perp_Q) \cong
O(\Hom_{\Stab_G(u)}(Q, (T_u U)^\perp)),
\]
 the orthogonal group acting
on the associated vector space $\Hom_{\Stab_G(u)}(Q, (T_u U)^\perp))$,
i.e., the multiplicity space of $(T_u U)^\perp_Q$. As explained in
\emph{op.~cit.}, the composition $\pi_1(X) \to \prod_{Q \in R_q((T_u
  U)^\perp)} \ZZ/2 \onto \ZZ/2$ with the projection to the factor $Q$
is nothing but application of the first Stiefel-Whitney class
$w_1(\Hom_{\Stab_G(u)}(Q, (T_u U)^\perp))$ of this orthogonal vector
bundle.

Now, let $\widetilde Y$ be the cover of $Y$ corresponding to the
kernel of the map $\pi_1(Y) \to H_1(Y, \ZZ/2)$.  Set $U := (\widetilde
Y \times V)/H_1(Y, \ZZ/2)$, where $H_1(Y, \ZZ/2)$ acts as follows.
First, it acts by the defining action on the factor of $\widetilde
Y$.  Next, for each $Q \in R_q(V)$, fix an isomorphism $V_Q \cong
Q^{m_Q}$. Then, let each $\gamma \in H_1(Y, \ZZ/2)$ act on $V_Q \cong
Q^{m_Q}$ by $\pm \Id \oplus \Id^{m_Q-1}$, where the sign is the image
of $\gamma$ under the composite map $\pi_1(Y) \to \prod_{Q' \in
  R_q(V)} \ZZ/2 \onto \ZZ/2 \cong \{\pm 1\}$ corresponding to $Q$.
Taking the direct sum, we obtain an action of $H_1(Y, \ZZ/2)$ on $V$,
and taking the product with the defining action on $\widetilde Y$,
we obtain an action of $H_1(Y, \ZZ/2)$ on $\widetilde Y \times V$.

It follows from the construction that $X_0 := Y$ is a symplectic leaf
of $X := U/G$: since $V^G = \{0\}$, $X_0 = U^G$.  Moreover, for $x \in
X_0$ and $u \in U$ mapping to $x$, $(T_u U)^\perp = U$.  It is
straightforward to check that the resulting map (b) is the given one.
\end{proof}

\subsection{Acknowledgements}  We are grateful to Victor Ginzburg
for useful discussions.
The first author's work was
partially supported by the NSF grant DMS-1000113. The second author is
a five-year fellow of the American Institute of Mathematics, and was
partially supported by the ARRA-funded NSF grant DMS-0900233.

\section{Proofs}
\label{typeas}
\subsection{Proofs of theorems and corollaries from \S
  \ref{ss:typea-int-dmod}}
\subsubsection{Proof of Theorem
  \ref{sympdmthm}} \label{ss:sympdmthm-pf} The symplectic leaves of $X
:= S^n Y$ are exactly the quotients $\{z \in Y^n \mid \Stab_{S_n}(z) =
G\} / N_{S_n}G$, where $G < S_n$ has the form $G = S_{i_1}^{r_1}
\times \cdots \times S_{i_k}^{r_k}$, and $N_{S_n}G < S_n$ is its
normalizer. Here, we may choose $i_1 < i_2 < \cdots < i_k$.  By
\cite[Theorem 3.1]{ESdm} and its proof, $M(X)$ is a holonomic
$\caD$-module whose composition factors are intermediate extensions of
local systems (i.e., $\cO$-coherent $\caD$-modules) on the symplectic
leaves.  Let $G = S_{i_1}^{r_1} \times \cdots \times S_{i_k}^{r_k}$ be
a fixed subgroup of $S_n$ as above (with $r_1 i_1 + \cdots + r_k i_k =
n$).  Let $(X^G)^{\circ} \subseteq X$ be the corresponding symplectic
leaf, and let $X^G$ denote its closure. Set $U := X \setminus (X^G
\setminus (X^G)^{\circ})$.  One has an obvious surjection $M(U) \onto
\Omega_{(X^G)^{\circ}}$ sending $1$ to the volume form.  As a result,
the intermediate extension of $\Omega_{(X^G)^{\circ}}$ is a
composition factor of $M(X)$.  To deduce the desired result,
therefore, it suffices to show that these are all the composition
factors, occurring with multiplicity one, and that $M(X)$ is
semisimple.

In order to prove that $M(X)$ is semisimple, we prove a more general
result: let $\phi: Y^n \onto X$ be the defining surjection, and
consider the $\caD$-module $M_\phi(Y^n)$ defined in \cite{ESdm}: this
is the quotient of $D_{Y^n}$ by the right ideal generated by
Hamiltonian vector fields of the form $\xi_{\phi^* f}$ for $f \in
\cO_X = \cO_Y^{S_n}$.  According to \cite[Theorem 3.1]{ESdm},
$M_\phi(Y^n)$ is a holonomic $\caD$-module on $Y^n$, and moreover
$(\pi_* M_\phi(Y^n))^{S_n} \cong M(X)$.  Thus, it suffices to show
that $M_\phi(Y^n)$ is semisimple.

To prove this, we recall again from \cite[Theorem 3.1]{ESdm} and its
proof that the singular support of $M_\phi(Y^n)$ in $T^* (Y^n)$ is
contained in the locus of pairs $(z, v)$ with $z \in Y^n, v \in T^*_z
Y^n$, such that $v \cdot \xi_{\phi^* f}|_{v} = 0$ for all $f \in
\cO_{X}$.  This is the union of the conormal bundles of the inverse
images of symplectic leaves on $X$.  Specifically, the closure of the
inverse image of each symplectic leaf is of the form $Y^{r_1 + \cdots
  + r_k} \subseteq (Y^{i_1})^{r_1} \times \cdots \times
(Y^{i_k})^{r_k} = Y^n$. Hence, the composition factors of
$M_\phi(Y^n)$ are $S_n$-equivariant local systems on these smooth,
closed subvarieties.  

We claim that $\Ext^1$ between any two such $\caD$-modules supported
on distinct diagonals is trivial.  Since the singular supports of
these $\caD$-modules are the conormal bundles of the given smooth
symplectic subvarieties, the claim follows from the more general
\begin{lemma} \label{l:extsubvars} Suppose that $Z$ is a smooth
  variety, and $Z_1, Z_2 \subseteq Z$ as well as $Z_1 \cap Z_2$ are
  smooth closed subvarieties, all of pure dimension.  Let
  $\mathcal{L}_1, \mathcal{L}_2$ be local systems on $Z_1$ and $Z_2$,
  respectively, and let $i_1: Z_1 \to Z$ and $i_2: Z_2 \to Z$ be the
  inclusions.  Then,
\begin{equation}\label{e:extsubvars}
  \Ext^j((i_1)_* \mathcal{L}_1, (i_2)_* \mathcal{L}_2) = 0, \text{ for } j < 
  (\dim Z_1 - \dim Z_1 \cap Z_2) + (\dim Z_2 - \dim Z_1 \cap Z_2).
\end{equation}
\end{lemma}
Namely, the result follows from the lemma since, in our case, $Z_1$,
$Z_2$, and $Z_1 \cap Z_2$ are all even dimensional and $Z_1 \neq Z_2$.
\begin{proof}[Proof of Lemma \ref{l:extsubvars}]
By adjunction, the LHS of \eqref{e:extsubvars} identifies with 
\begin{equation}\label{e:ext1}
\Ext^j(i_2^* (i_1)_* \mathcal{L}_1, \mathcal{L}_2).
\end{equation}
Next, let $i_{12,k}: Z_1 \cap Z_2 \to Z_k$ be the inclusions for $k \in \{1,2\}$.  Then, applying proper base change for the closed
embedding $i_{12,2}$, we can rewrite  \eqref{e:ext1} as
\begin{equation}
\Ext^{j}((i_{12,2})_* i_{12,1}^* \mathcal{L}_1, \mathcal{L}_2).
\end{equation}
Since $i_{12,2}$ is a closed embedding, $(i_{12,2})_*
=(i_{12,2})_!$. Applying adjunction, we obtain
\begin{equation}
\Ext^{j}(i_{12,1}^* \mathcal{L}_1, i_{12,2}^! \mathcal{L}_2).
\end{equation}
Now, $i_{12,1}^* \mathcal{L}_1$ is a local system shifted by $-(\dim Z_1 - \dim Z_1 \cap Z_2)$, and $i_{12,2}^! \mathcal{L}_2$ is a local system shifted by $\dim Z_2 - \dim Z_1 \cap Z_2$. 
 So, the above vanishes when $j < (\dim Z_1 - \dim Z_1 \cap Z_2) + (\dim Z_2 - \dim Z_1 \cap Z_2)$ (or when $j >  (\dim Z_1 - \dim Z_1 \cap Z_2) + (\dim Z_2 - \dim Z_1 \cap Z_2) + \dim Z_1 \cap Z_2$).  
\end{proof}
\begin{remark}
  In fact, the above lemma is needed for the omitted proof of
  \cite[Theorem 4.21]{ESdm}.  So, even though we could have deduced
  semisimplicity from that theorem, the above argument cannot be
  avoided.
\end{remark}
It remains to prove that the intermediate extensions
$\Omega_{(X^G)^\circ}$ are all of the composition factors of $M(X)$,
and that they occur with multiplicity one.  Then, the irreducible
composition factors of $M_\phi(Y^n)$ are all supported on distinct
diagonal subvarieties of $Y^n$, so the above argument implies that
$M_\phi(Y^n)$, and hence $M(X)$, are semisimple.  Since the
composition factors are exactly the claimed direct summands of $M(X)$,
the theorem also follows.

So, we prove that the intermediate extensions of
$\Omega_{(X^G)^\circ}$, i.e., the IC $D$-modules of $(X^G)^\circ$, are
all of the composition factors of $M(X)$, and that they occur with
multiplicity one.  It suffices to consider the formal neighborhood of
a point of $(X^G)^\circ$.  Then, the computation reduces to the case
that $G = S_n$ and $X^G = (X^G)^\circ = Y \subseteq S^n Y$, and
moreover, we may reduce to the case that $Y = V$ is a symplectic
vector space, and consider the formal neighborhood of zero,
$\widehat{\cO_V} = \CC\dblsqbrs{x_1, \ldots,x_d,y_1,\ldots, y_d}$.
Let $\delta_V$ be the delta-function $\caD$-module of the diagonal $V
\subseteq S^n V$.  Since $V$ is now a symplectic vector space
\cite[Corollary 4.16]{ESdm} implies that $M(S^n V)$ is semisimple, and
a direct sum of IC $D$-modules of the symplectic leaves with some
multiplicities (in fact, \emph{op.~cit.} implies that the multiplicity
of $\delta_V$ is $\dim \HP_0(\cO_{V^{n-1}}^{S_n})$, which would reduce
us to Theorem \ref{snp1thm}, but we will instead deduce that theorem
from the present one).  It suffices to prove
\begin{equation}
\Hom_{D_{S^n V}}(M(S^n V), \delta_{V}) \cong
\CC.
\end{equation}
This may be restated and proved without the use of $\caD$-modules:
\begin{lemma} \label{l:hidist} The space of symmetric polydifferential
  operators $\psi: (\mathcal{O}_V)^{\otimes (n-1)} \rightarrow \cO_V$
  invariant under Hamiltonian flow is one-dimensional, and spanned by
  the multiplication map.
\end{lemma}
Note that, actually, we only need to show that there are no
$S_n$-invariant operators, with the $S_n$ action given by viewing the
polydifferential operators in the lemma as distributions on $n$
functions; the lemma is a slightly more general result, requiring only
$S_{n-1}$-invariance.

We remark that this lemma is tantamount to Theorem \ref{snp1thm},
i.e., one can directly show that the above space of polydifferential
operators is identified with $\HP_0(\cO_{V^{n-1}/S_n}, \cO_{V})$ (at least
if we require that the operators be $S_n$-invariant).  For
details, see \cite[\S 4]{ReScmat}.

We further remark that the space mentioned in the lemma can
alternatively be viewed as the space of $C^\infty$
Hamiltonian-invariant distributions on $S^n V$ supported on the
diagonal, since finite-dimensionality guarantees that $\Hom_{D_{S^n
    V}}(M(S^n V), \delta_{V})$ is the same when considered in the
$C^\infty$ context.  Then, a polydifferential operator $\psi$ of
degree $n-1$ becomes a distribution $\Psi$ on $V^n$ by the
prescription $\Psi(f_1, \ldots, f_n) = \int \psi(f_1, \ldots, f_{n-1})
f_n$. They are supported on the diagonal since they depend only on
(finitely many) partial derivatives of $f_1 \times \cdots \times f_n$
evaluated at the diagonal.

\begin{proof}
  It suffices to pass to the formal completion and consider
  polydifferential operators on $\widehat{\cO_V}$.  Such
  polydifferential operators are determined by their value on elements
  $f^{\otimes (n-1)}$ for $f \in \widehat{\cO_V}$, since they are
  symmetric and hence determined by their restriction to
  $\Sym^{n-1} \widehat{\cO_V}$. Furthermore, we can assume that $f'(0)
  \neq 0$, since the complement of this locus in the pro-vector space
  $\widehat{\cO_V}$ has codimension equal to $\dim V \geq 2$.

  Next, by the formal Darboux theorem, by applying a formal
  symplectomorphism of $V$, we may assume $f = x_1$.  Since all formal
  symplectomorphisms are obtained by integrating Hamiltonian vector
  fields, it suffices to consider the value $\psi(x_1^{\otimes
    (n-1)})$.  This value must be a function that depends only on
  $x_1$, since these are the only functions invariant under all
  symplectomorphisms fixing $x_1$. By linearity and invariance under
  conjugation by rescaling $x_1$ (and applying the inverse scaling to
  $y_1$), we deduce that $\psi(x_1^{\otimes (n-1)}) = \lambda \cdot
  x_1^{n-1}$ for some $\lambda \in \CC$.  Thus, on $x_1^{\otimes (n-1)}$, 
$\psi$ coincides with
  $\lambda$ times the multiplication operator, $f_1 \otimes \cdots \otimes f_{n-1} \mapsto \lambda f_1 \cdots f_{n-1}$.  The latter operator is evidently
symmetric
 and invariant under Hamiltonian flow. On the other hand, we have argued that a symmetric operator invariant under Hamiltonian flow is uniquely determined by its value on $x_1^{\otimes (n-1)}$.  So $\psi$ is equal to $\lambda$ times the
multiplication operator, as desired.
\end{proof}
\subsubsection{Proofs of Corollaries \ref{symphp0cor} and \ref{c:defcor1-dr}}
\begin{proof}[Proof of Corollary \ref{symphp0cor}]
By definition, for all (affine) Poisson varieties $X$, 
$\HP_\bullet^{DR}(X) = \pi_* M(X)$, where $\pi_*$ is the \emph{derived} pushforward of the $\caD$-module $M(X)$ to a point. This identifies the first term in \eqref{hpdreqn1} with the LHS of \eqref{e:sympdmthm}.

Next, since $Y$ is symplectic, by \cite[Example 2.6]{ESdm}, $M(Y) =
\Omega_Y$.  Therefore, $\pi_*((\Delta_{i})_* (\Omega_Y)^{\boxtimes r})
\cong (\HP_{\bullet}^{DR}(Y))^{\otimes r}$, with the canonical $S_r$
action given by permutation of components.  This identifies the second
term in  \eqref{hpdreqn1} with the RHS of
\eqref{e:sympdmthm}.

It remains to consider the third term in \eqref{hpdreqn1}.  Here, we
use again the fact that $M(Y) = \Omega_Y$, together with the standard
fact that, since $Y$ is smooth and connected, $H_\bullet(\pi_*
\Omega_Y) \cong H^{\dim Y - \bullet}(Y)$.
\end{proof}
\begin{proof}[Proof of Corollary \ref{c:defcor1-dr}]
It suffices to show that the third term in \eqref{hpdreqn1} is identified
with the first term in \eqref{e:defcor1-dr}.  By results of
  Nest-Tsygan \cite{NeTs}, 
one has an isomorphism 
\begin{equation}\label{e:nets}
H^{\dim Y - \bullet}(Y)\dblparens{\hbar} \cong
 \HP_{\bullet}(Y)\dblparens{\hbar} \iso 
\HH_{\bullet}(A_\hbar[\hbar^{-1}]),
\end{equation}
where $\HP_\bullet(Y)$ is the usual Poisson homology (which is well
known to be isomorphic to $H^{\dim Y - \bullet}(Y)$ when $Y$ is
symplectic, since the Poisson homology complex identifies with the de
Rham complex). Thus, it suffices to show that
\begin{equation}
\bigoplus_{n \geq 0} \HH_\bullet(\Sym^n A_\hbar)[\hbar^{-1}]^* \cong \Sym(\HH_\bullet(A_\hbar[\hbar^{-1}])^*[t]),
\end{equation}
again taking the $\CC\dblparens{\hbar}$-linear dual.  Since $A_\hbar$ is an
infinite-dimensional, simple algebra with trivial center, this follows
from \cite[Corollary 3.3]{EO}. Since we will need this again later, we
state it below.
\end{proof}
We used here and will continue to use the following result from
\cite{EO}, which we state somewhat more explicitly than is in
\emph{op.~cit.} (we omit the proof of the more explicit formula, as we
do not essentially need it):
\begin{theorem}\cite[Corollary 3.3]{EO} \label{eothm} Let $A$ be an
  infinite-dimensional simple algebra over a field of
characteristic zero with trivial center.
Then, the coalgebra $\HHH_\bullet(A) := \bigoplus_{n \geq 0} \HH_\bullet(\Sym^n A)$ is a polynomial coalgebra,
\begin{equation} \label{eofla}
\Psi: \HHH_\bullet(A) \iso \Sym (\HH_\bullet(A)[t]),
\end{equation}
where the isomorphism is the unique coalgebra map which is graded with
respect to $|\HH_\bullet(A)|=|t|=1$  such that, for every $n$,
composition with 
the projection to $t^{n-1} \HH_\bullet(A)$ restricts on $\HH_\bullet(\Sym^n A)$ to a map
\begin{equation}
\HH_\bullet(\Sym^n A) \rightarrow \Sym (\HH_\bullet(A)[t]) \onto t^{n-1} \HH_\bullet(A)
\end{equation}
of the form, in Hochschild degree zero,
\begin{equation} \label{e:eofla-eqn}
[a_1 \& \cdots \& a_n]) \mapsto \frac{1}{n!}  \sum_{\sigma \in S_n}t^{n-1} [a_{\sigma(1)} \cdots a_{\sigma(n)}],
\end{equation}
and similarly is the natural multiplication map on Hochschild $m$-chains for all $m \geq 0$,
\[
c_0 \otimes \cdots \otimes c_m \mapsto c_0' \otimes \cdots
\otimes c_m',
\]
where $c_i \mapsto c_i'$ is the map \eqref{e:eofla-eqn}.
\end{theorem}

\subsection{Proofs of theorems and corollaries from \S \ref{ss:typea-int}}
\subsubsection{Proofs of Theorems \ref{symphp0thm} and \ref{snp1thm}}
Theorem \ref{symphp0thm} already follows from the corollary
\ref{symphp0cor} of Theorem \ref{sympdmthm}, so it remains only to
prove Theorem \ref{snp1thm}.
\begin{proof}[Proof of Theorem \ref{snp1thm}]
  As in the introduction, write $S^n V \iso (V \times V^{n-1}/S_n)$
  where the map to the first factor is given by averaging.  We deduce
  that $M(S^n V) \cong \Omega_V \boxtimes M(V^{n-1}/S_n)$.  Recall
  from \cite[Theorem 4.13]{ESdm} that, for any symplectic vector space
  $U$ and finite subgroup $G < \Sp(U)$, the space $\HP_0(\cO_{U}^G)$
  naturally identifies with the multiplicity space of the
  delta-function $\caD$-module of the origin in $M(U/G)$, which is
  semisimple.\footnote{For general affine Poisson varieties $X$ and $x
    \in X$, the space $\Hom_{D_X}(M(X), \delta_x)$ identifies with a
    subspace of $\HP_0(\cO_X)^*$; when $\cO_X$ is nonnegatively graded
    with the ideal of $x$ as the augmentation ideal, this is an
    equality.  
% Moreover, for $X=U/G$, $M(X)$ is semisimple.
  } Hence, it also identifies with the multiplicity space of the
  delta-function $\caD$-module of the diagonal $V \subset S^n V$ in
  $M(S^n V)$.  This multiplicity space is one-dimensional (in fact,
  the main step of the proof of Theorem \ref{sympdmthm} was to show
  this).
\end{proof}
In the appendix, we will give a different, elementary proof of Theorem
\ref{snp1thm}.  A proof without using $M(S^n V)$, requiring only the
Darboux theorem, can also be obtained from Lemma \ref{l:hidist}
following the comments after the statement of the lemma.

\subsubsection{Proofs of Corollaries \ref{c:mpfla}--\ref{defcor2}}
\begin{proof}[Proof of Corollary \ref{c:mpfla}]
  This is immediate by expanding the LHS of \eqref{e:symphp0thm},
  since $\HP_0(\cO_{Y}) \cong H^{\dim Y}(Y)$: namely, the subspace of
  the LHS spanned by terms of the form $f_1t^{r_1} \& \cdots \&
  f_mt^{r_m}$ for fixed $r_1 \leq r_2 \leq \cdots \leq r_m$ has a
  basis by requiring that the $f_i$ lie in a fixed basis of
  $\HP_0(\cO_Y)$; then the number of these is the number of
  $\dim \HP_0(\cO_Y)$-multipartitions $(\lambda_1, \ldots, \lambda_{\dim \HP_0(\cO_Y)})$ of $n$,
  i.e., $|\lambda_1| + \cdots + |\lambda_{\dim \HP_0(\cO_Y)}| = n$, such that there are a total of $m$ parts appearing in all the partitions, of lengths $r_1+1, \ldots, r_m+1$.  (In particular, $(r_1+1) + \cdots + (r_m+1) = n$.)
\end{proof}
\begin{proof}[Proof of Corollary \ref{defcor1}]
  This is a direct consequence of Corollary \ref{c:defcor1-dr} (or we
  can prove it in the same manner, using only Theorem \ref{symphp0thm}
  rather than Theorem \ref{sympdmthm}).
\end{proof}
\begin{proof}[Proof of Corollary \ref{defcor2}]
  This is an immediate consequence of Corollary \ref{c:mpfla}, using
  the canonical surjection $\HP_0(\cO_{S^n Y})\dblparens{\hbar} \onto
  \HH_0(B_\hbar[\hbar^{-1}])$.
\end{proof}
\subsubsection{Proofs of Corollaries \ref{defcor3}--\ref{defcor6}}
\begin{proof}[Proof of Corollary \ref{defcor3}]
  By Theorem \ref{snp1thm}, $\HP_0(\cO_{V^n}^{S_{n+1}}) \cong \CC$; it
  suffices to show that $\HH_0(\caD_{U^n}^{S_{n+1}}) \cong \CC$.  This
  is a consequence of \cite{AFLS}: the dimension of
  $\HH_0(\caD_{U^n}^{S_{n+1}})$ is equal to the number of conjugacy
  classes of elements in $S_{n+1}$ which act without eigenvalue one on
  $U^n$; there is exactly one such conjugacy class, namely the
  conjugacy class of the $(n+1)$-cycle.
\end{proof}
% \begin{proof}[Proof of Corollary \ref{defcor4}]
%   This is an immediate consequence of Theorem \ref{snp1thm} and the
%   canonical surjection $\HP_0(\cO_{V^n}^{S_{n+1}}) \onto \gr
%   \HH_0(B)$.
% \end{proof}
\begin{proof}[Proof of Corollary \ref{defcor5}]
This is Corollary \ref{defcor3} in the case that $\dim V = 2$.
\end{proof}
\begin{proof}[Proof of Corollary \ref{defcor6}]
  By upper semicontinuity of $\dim \HH_0(B)$ in the family of filtered
  quantizations $B$, it suffices to
  show that, for generic spherical rational Cherednik algebras $B$
  deforming $\cO_{\CC^{2n}}^{S_{n+1}}$, $\dim \HH_0(B) = 1$.  By
  \cite[Theorem 1.8]{EGsra}, for generic $B$, $\HH_\bullet(B) \cong
  \HH_\bullet(\caD_{\CC^n}^{S_{n+1}})$. Thus, the result follows from
  \cite{AFLS}, as already explained in the proof of Corollary
  \ref{defcor3} (or one can simply refer to that corollary or
  Corollary \ref{defcor5}).
\end{proof}

\subsubsection{Proofs of Corollaries \ref{c:los1} and \ref{c:los2}}
\begin{proof}[Proof of Corollary \ref{c:los1}]
  Losev's \cite[Appendix A]{ESdm} implies the following result.  Let
  $X$ be an affine Poisson variety with finitely many (locally closed)
  symplectic leaves $X_1, \ldots, X_k$.  Let $B_\hbar$ or $B$ be a
  deformation or filtered quantization of $\cO_X$ (the latter only in
  the case that $\cO_X$ is nonnegatively graded).  For each symplectic
  leaf $X_i$ let $x_i \in X_i$ be a point. Let $\hat \cO_{X,x_i}$ be
  the formal completion of $\cO_X$ at $x_i$.  Now, write $\hat X_{x_i}
  := \Spf \hat \cO_{X, x_i}$ for the formal neighborhood of $x_i$ in
  $X$, where $\Spf$ refers to the ``formal'' spectrum of prime ideals
  in $\cO_{X, x_i}$ which are closed under the
  $\mathfrak{m}_{x_i}$-adic topology, and $\mathfrak{m}_{x_i}$ is the
  maximal ideal associated to $x_i \in X$. According to
  \cite[Proposition 3.3]{Kalss}, there is an isomorphism $\hat X_{x_i}
  \cong \hat (X_i)_{x_i} \hat \times \hat Z_i$, for some ``slice''
  subvariety $\hat Z_i \subseteq \hat X_{x_i}$.  That is, $\hat
  \cO_{X, x_i} \cong \hat \cO_{X_i,x_i} \hat \otimes \cO_{\hat Z_i}$,
  where $\cO_{\hat Z_i}$ is a quotient of $\cO_{X,x_i}$ by a complete
  ideal, and $\hat \otimes$ denotes the completed tensor product.

  We will need to consider the space $\HP_0(\cO_{\hat Z_i})=\cO_{\hat
    Z_i} / \{\cO_{\hat Z_i}, \cO_{\hat Z_i}\}$. (Note that, as pointed
  out in \cite[Proposition 3.10]{ESdm}, $\{\cO_{\hat Z_i}, \cO_{\hat
    Z_i}\}$ is a closed subspace of $\cO_{\hat Z_i}$ in the adic
  topology, and $\HP_0(\cO_{\hat Z_i}) = \lim_{n \to \infty}
  \HP_0(\cO_{\hat Z_i} / \mathfrak{m}_{x_i}^n)$.)

  The following is a direct consequence of \cite[Appendix A]{ESdm}, using
  \cite[Proof of Corollary 3.13]{ESdm}:
\begin{theorem}\cite[Appendix A]{ESdm} 
Every prime ideal of $B_\hbar[\hbar^{-1}]$ over $\CC\dblparens{\hbar}$
is supported on $\overline{X_i}$ for some $i$. For each $i$, the number
of such ideals is at most $\dim_{\CC\dblparens{\hbar}} \HP_0(\cO_{\hat Z_i})\dblparens{\hbar}$.
\end{theorem}
Now, the closures of the symplectic leaves of $S^n Y$ are exactly the
images of all possible compositions $Y^m \to (S^{i_1} Y)^{r_1} \times
\cdots \times (S^{i_k} Y)^{r_k} \to S^n Y$, where $m = r_1 + \cdots +
r_k$, $n = r_1 i_1 + \cdots + r_k i_k$, the first map is the product
of the $m$ diagonal embeddings $\Delta_{i_1}^{r_1} \times \cdots
\times \Delta_{i_k}^{r_k}$, and the second map is the obvious
projection. At a point $x_i$ of the locally closed symplectic leaf
$X_i$ with this closure, the slice $\hat Z_i$ such that $(\widehat{S^n
  Y})_{x_i} \cong \hat X_{x_i} \hat \times \hat Z_i$ can be taken to
be isomorphic to the formal neighborhood of the origin in
$(V^{i_1-1}/S_{i_1})^{r_1} \times \cdots \times
(V^{i_k-1}/S_{i_k})^{r_k}$, where $V$ is a symplectic vector space of
dimension equal to $\dim Y$. Namely, if $\phi: Y^n \to S^n Y$ is the
projection, we can consider a preimage $\widetilde x_i \in S^n Y$ of
$Y^n$ and look at the completed conormal fiber of $\phi^{-1}(X_i)$ at
$\widetilde x_i$, then project back down to $S^n Y$, to get $\hat
Z_i$.

Therefore, $\HP_0(\cO_{\hat Z_i})
\cong \HP_0(V^{i_1-1}/S_{i-1})^{\otimes r_1} \otimes \cdots \otimes \HP_0(V^{i_k-1}/S_{i-1})^{\otimes r_k} \cong \CC$.  So, there is at most one prime ideal supported on $\overline{X_i}$.  Note that $\operatorname{codim}_{S^n Y}(\overline{X_i}) = (n-m)\dim Y$, where $m = r_1 + \cdots + r_k$ as above.

Thus, the number of prime ideals with support of codimension $(n-m) \dim Y$
is at most the number of partitions of $n$ with $m$ parts. This immediately implies the statement.
\end{proof}
\begin{proof}[Proof of Corollary \ref{c:los2}]
  Note that $V^{n+1}/S_{n+1} \cong V \times V^n/S_{n+1}$, with the
  projection $V^{n+1}/S_{n+1} \to V$ given by averaging the $n+1$
  elements of $V$ in the ordered $(n+1)$-tuple, and the map
  $V^{n+1}/S_{n+1} \to V^n/S_{n+1}$ given by subtracting the average
  from each element of the $(n+1)$-tuple.  Therefore, the symplectic
  leaves of $V^{n+1}/S_{n+1}$ are all of the form $V \times X_i$ where
  $X_i$ is a symplectic leaf of $V^n/S_{n+1}$, and this establishes a
  bijection between the symplectic leaves of $V^{n+1}/S_{n+1}$ and
  those of $V^n/S_{n+1}$.  The corollary then follows from Corollary
  \ref{c:los1}.
\end{proof}

\subsubsection{Proofs of Corollaries \ref{c:fpoiss1}--\ref{c:fpoiss2}}
\begin{proof}[Proof of Corollary \ref{c:fpoiss1}]
  This follows from Corollary \ref{c:mpfla} and the
  surjection $\HP_0(\cO_{S^n Y})\dblparens{\hbar} \onto \gr
  \HP_0(B_\hbar[\hbar^{-1}])$.
\end{proof}
\begin{proof}[Proof of Corollary \ref{c:fpoiss2}]
  By the comments before the corollary, it suffices to show that
  $\dim \HP_0(B) = 1$ when $B$ is a commutative spherical Cherednik
  algebra which is a filtered Poisson deformation of
  $\cO_{\CC^{2n}}^{S_{n+1}}$. The fact that $\dim \HP_0(B) \leq 1$
  follows from the surjection $\HP_0(\cO_{V^n}^{S_{n+1}}) \onto \gr
  \HP_0(B)$ and Theorem \ref{snp1thm}. For the opposite inequality, by
  upper semicontinuity of $\dim \HP_0(B)$, it suffices to show that,
  for generic commutative spherical Cherednik algebras $B$ deforming
  $\cO_{\CC^{2n}}^{S_{n+1}}$, one has $\HP_0(B) \neq 0$. This result
  follows because, by \cite[Corollary 1.14]{EGsra}, generic spherical
  rational Cherednik algebras $B$ deforming $\cO_{\CC^{2n}}^{S_{n+1}}$
  are of the form $B = \cO_X$ for $X$ smooth and symplectic, with
  one-dimensional top cohomology, i.e., $\HP_0(B) \cong
  \CC$. (Alternatively, without using that $\Spec B$ is smooth for
  generic $B$, one could take a deformation quantization $B_\hbar$ of
  $B$ such that $B_{\hbar}[\hbar^{-1}]$ is isomorphic to a
  noncommutative spherical Cherednik algebra over $\CC\dblparens{\hbar}$
  deforming $\cO_{\CC^{2n}}^{S_{n+1}}\dblparens{\hbar}$, so $\dim
  \HH_0(B_{\hbar}[\hbar^{-1}]) = 1$ by \cite[Theorem 1.8]{EGsra} (or
  by Corollary \ref{defcor6}, which uses \emph{op.~cit.}).  Then, one
  concludes using the canonical surjection $\HP_0(B)\dblparens{\hbar} \to \gr
  \HH_0(B_{\hbar}[\hbar^{-1}])$.)
\end{proof}

\appendix

\section{Type $D$ Weyl groups, by T. Schedler}\label{s:typed-int}
In this appendix, we compute $\HP_0(\cO_{\CC^{2n}}^{D_n})$, where $D_n
< \GL(\CC^n) < \Sp(\CC^{2n})$ is the type $D_n$ Weyl subgroup. Recall
that $D_n = S_n \ltimes (\ZZ/2)^{n-1}$, and we let $\CC^n$ be its
reflection representation, where $S_n$ acts by permuting components,
and $(\ZZ/2)^{n-1}$ acts by diagonal matrices whose diagonal entries
are $\pm 1$ which have determinant one (i.e., an even number of $-1$
entries).

Note that $D_n$ is an index-two subgroup of $B_n = C_n = S_n \ltimes
(\ZZ/2)^n$.  Also, $\bigoplus_{n \geq 0}
\HP_0(\cO_{\CC^{2n}}^{B_n})^*$ is a bigraded algebra, graded by the
symmetric power degree, $n$, and the weight degree (degree of
polynomials in $\cO_{\CC^{2n}}$ for all $n$).  Recall from
\cite{ESsym}:
\begin{theorem}\label{bnhp0thm}\cite{ESsym}\footnote{Note that \cite[Theorem 1.1.3]{ESsym}
    is for the much more general situation of symmetric powers of
    isolated surface singularities in $\CC^3$ with a contracting
    $\CC^*$-action, but we only need the case of the surface
    $\CC^2/(\ZZ/2)$.} There is an isomorphism of bigraded algebras
\begin{equation}\label{bnhp0eqn}
\bigoplus_{n \geq 0} \HP_0(\cO_{\CC^{2n}}^{B_n})^* \cong \CC[s_1, s_2, \ldots],
\end{equation}
where $s_i$ has symmetric power degree $i$ and weight $4(1-i)$.\footnote{We assign $s_i$ \emph{nonpositive} weight because it lies in the dual space
to $\HP_0(\cO_{\CC^{2n}}^{B_{i+1}})$, which is assigned nonnegative weight.}
\end{theorem}
Here, the algebra structure on the LHS arises from the symmetrization
map: precisely, given $\phi \in (\Sym^m \cO^{\ZZ/2}_{\CC^2})^* =
(\cO_{\CC^{2m}}^{B_m})^*$ and $\psi \in (\Sym^n \cO^{\ZZ/2}_{\CC^2})^* =
(\cO_{\CC^{2n}}^{B_n})^*$, then $\phi \cdot \psi \in (\Sym^{m+n}
\cO_{\CC^2}^{\ZZ/2})^*$ is defined by $\phi \cdot \psi = \phi \boxtimes \psi$,
viewing $\Sym^{m+n} \cO^{\ZZ/2}_{\CC^2}$ as the subspace of $T^{m+n}
\cO_{\CC^2}^{\ZZ/2}$ of symmetric tensors.

We now compute $\HP_0(\cO_{\CC^{2n}}^{D_n})$.  Let us view $s_i$ as
the coordinate functions on the infinite-dimensional space
$\CC\dblsqbrs{x^2}$ (we will explain why in the proof), so that, for all $f
\in \CC\dblsqbrs{x^2}$,
\[
f = s_1(f) + s_2(f) x^2 + s_3(f) x^4 + \cdots.
\]
 We need to define certain vector
fields $\xi_k$ on $\CC\dblsqbrs{x^2}$ for $k \geq 1$.  First, let $Q(z)$ be the Taylor series of $\sqrt{1+z}$, i.e., $Q(z) = 1 + \frac{z}{2} - \frac{z^2}{8} + \frac{z^3}{16} - \cdots$. Then, we define
\begin{equation}
\xi_k := \mathcal{V}\biggl(\frac{d}{dx} \Bigl(x^{2k-1} \cdot Q\bigl(\frac{1}{s_{2k}}\sum_{i \geq 2k+1}
s_i x^{2(i-2k)}\bigr) \Bigr)\biggr),
\end{equation}
where \[
\mathcal{V}\bigl(\sum_{i \geq 0} f_i x^{2i}\bigr) := \sum_{i \geq 0} f_i \partial_{s_{i+1}}, \quad f_i \in \CC[s_1, s_2, \ldots].
\]
Explicitly, the first few terms of $\xi_k$ can be written out as
\begin{equation}
\xi_k = \bigl( (2k-1) \partial_{s_{k}} +  \frac{2k+1}{2} \frac{s_{2k+1}}{s_{2k}} \partial_{s_{k+1}} + \cdots \bigr),
\end{equation}
where here $\cdots$ means terms that are multiples of $s_{2k+j}$ for
$j \geq 2$.
\begin{theorem}\label{dnhp0thm}
The sum $\bigoplus \HP_0(\cO_{\CC^{2n}}^{D_n})^*$ is naturally a bigraded subalgebra of $\bigoplus_{n \geq 0} \HP_0(\cO_{\CC^{2n}}^{B_n})^*$. In terms of \eqref{bnhp0eqn}, it is the subalgebra of elements $f$ such that, for all $k \geq 1$, 
\begin{equation} \label{e:xik-cond}
\xi_k(f)|_{s_1 = \cdots = s_{2k-1}=0, s_{2k} \neq 0} = 0.
\end{equation}
\end{theorem}
\begin{remark}\footnote{Thanks to P. Etingof for pointing out this
    observation.}
  It is interesting to try to integrate the above vector fields, in
  order to interpret solutions $f \in \bigoplus_{n \geq 0}
  \HP_0(\cO_{\CC^{2n}}^{D_n})^*$ as functions on $\CC[\![x^2]\!]$
  invariant under a certain flow. We can interpret this flow as
  follows: A curve $h(t)$ in $\CC[\![x^2]\!]$ is invariant if and only if $h_t =
  - (\sqrt{h})_x$, i.e., setting $u := 2 \sqrt{h}$, we should have
\[
u_x + u u_t = 0.
\]
This equation is the well known inviscid Burgers equation (with $t$
and $x$ swapped).  Then, the solutions should look like $u = f(t-ux)$
for some function $f$.

At $t=0$, we obtain $u(0,x) = f(-u(x,0) x)$.  So, in the case that
$u(x,0) \in x \CC[\![x^2]\!]^\times$, i.e., $h(0) \in \CC^\times \cdot
x^2 + x^4 \CC[\![x^2]\!]$, this implies that $f = \sqrt{g}$ where $g$
has linear behavior near $0$.  This implies that $u^2=g(t-ux)$, and
letting $G$ be an inverse of $g$, we can write $G(u^2) + ux - t = 0$,
which can now be solved for $u$.  For example, if $g(z)=-z$, then we
obtain $u(x,t) = \frac{x + \sqrt{x^2-4t}}{2}$.

However, it is not clear whether one can use this to simplify the
description of the algebra $\bigoplus_{n \geq 0}
\HP_0(\cO_{\CC^{2n}}^{D_n})^*$.
\end{remark}
\begin{corollary}
For $n \geq 7$, $\gr \HH_0(\caD_{\CC^{n}}^{D_n})^*$ is a proper
subspace of $\HP_0(\cO_{\CC^{2n}}^{D_n})^*$.  For $n \leq 6$, the
natural inclusion is an equality,
$\gr \HH_0(\caD_{\CC^{n}}^{D_n})^* = \HP_0(\cO_{\CC^{2n}}^{D_n})^*$.
\end{corollary}
\begin{proof} 
  Clearly it suffices to show that $\dim \HH_0(\caD_{\CC^{n}}^{D_n}) =
  \dim \HP_0(\cO_{\CC^{2n}}^{D_n})$ if and only if $n \leq 6$.  By the
  main result of \cite{AFLS}, for an arbitrary symplectic vector space
  $V$, $G < \Sp(V)$, and Lagrangian $U \subseteq V$, the dimension of
  $\HH_0(\caD_U^G)$ is equal to the number of conjugacy classes of
  elements $g \in G$ such that $g - Id$ is invertible (acting on
  $V$). (However, this says nothing about the \emph{filtration} on
  $\HH_0(\caD_U^G)$, which we deduce in this corollary.)  In the case
  at hand with $G = D_n$, the dimension of $\HH_0(\caD_{\CC^n}^{D_n})$
  therefore equals the number of partitions of $n$ with an even number
  of parts.

  Note that solutions of \eqref{e:xik-cond}, in particular, include
  all multiples of $s_1^2$.  One can inductively prove that, for $n >
  10$, there are more of the latter type of partitions than there are
  of the former.  Alternatively, more linearly independent solutions
  of \eqref{e:xik-cond} are given, for every monomial $g$ in $s_2,
  s_3, \ldots, s_{k+1}$, by $s_2^k s_{k+1} \cdot g - s_1 \xi_1(s_2^k
  s_{k+1} \cdot g)$ (this is a polynomial, and not merely a Laurent
  polynomial, because of the restriction on $g$).  One can inductively
  prove that the number of these plus the number of monomial multiples
  of $s_1^2$ exceed the number of even partitions of $n$ for $n > 8$
  and $n=7$; then it remains only to consider the case $n=8$, where
  one can find an additional solution not spanned by these (as
  reported in Figure \ref{f:expdn1}); it lies in weight -20.  The fact
  that the isomorphism stated in the corollary holds for $n \leq 6$ is
  a consequence of a straightforward explicit computation, or see
  Figure \ref{f:expdn1}.
\end{proof}

The above theorem, along with Theorem \ref{snp1thm} (for the type
$A_n$ cases) and the results of \cite{ESsym} (which imply the
$B_n=C_n$ cases) complete the computation of Poisson traces for
varieties $(\mathfrak{h} \oplus \mathfrak{h}^*)/W$ for $W$ one of the
classical series ($A, B=C$, and $D$) of finite Weyl groups and
$\mathfrak{h}$ its reflection representation.  Little is known about
the exceptional cases: only the case $G_2$ was computed in
\cite{AlFo}.  We also remark that, if we consider also the finite
Coxeter groups, the additional rank $\leq 3$ cases ($I_2(m)$ and
$H_3$) are computed in \cite{hp0bounds}. In all of these cases, one
has $\HP_0(\cO_{\mathfrak{h} \oplus \mathfrak{h}^*}^W) \cong \gr
\HH_0(\caD_{\mathfrak{h}}^W)$.

\subsection{Filtered quantizations and Poisson deformations}
Here we explain the analogous corollaries to those in the main body of
the paper, now for type $D_n$ rather than type $A_n$ Weyl groups.  For
all $n$, let $d_n$ be the dimension of $\HP_0(\cO_{\CC^{2n}}^{D_n})$,
as follows from the theorem (for $n \leq 34$, this can also be
obtained by evaluating the polynomials in Figures \ref{f:expdn1} and
\ref{f:expdn2} at $t=1$).  
The next corollary is an analogue of Corollary \ref{defcor3},
and is proved in the same manner:
\begin{corollary}
  Let $B$ be a filtered quantization of $\cO_{\CC^{2n}}^{D_n}$. Then, $\dim
  \HH_0(B) \leq d_n$, and the number of irreducible finite-dimensional
  representations of $B$ is at most $d_n$.
\end{corollary}
In particular, this includes the noncommutative spherical Cherednik
algebras deforming $\cO_{\CC^{2n}}^{D_n}$. (Note that we cannot obtain an
equality in this case since $\dim \HP_0(\cO_{\CC^{2n}}^{D_n}) > \dim
\HH_0(\caD_{\CC^n}^{D_n}) =$ the dimension of $\HH_0(B)$ for generic
noncommutative spherical Cherednik algebras deforming
$\cO_{\CC^{2n}}^{D_n}$. This is partly a reflection of the fact that $\CC^{2n}/D_n$
does not admit a symplectic resolution; see \S \ref{ss:csr} of the
main text.)

The next corollary is an analogue of Corollary \ref{c:fpoiss2}, proved
in the same manner:
\begin{corollary}
  Let $B$ be a filtered Poisson deformation of $\cO_{\CC^{2n}}^{D_n}$. Then,
  $\dim \HP_0(B) \leq d_n$, and the number of zero-dimensional
  symplectic leaves of $\Spec B$ is at most $d_n$.
\end{corollary}
In particular, this includes the commutative spherical Cherednik
algebras deforming $\cO_{\CC^{2n}}^{D_n}$. (For the same reason as before, we
cannot obtain an equality in this case.)

One can also formulate an analogue of Corollary \ref{c:los2} (which
can also be proved in the same manner; see also the proof of Theorem
\ref{dndmthm} in \S \ref{ss:dndmthm-proof}).  Recall the definition of
$p_{n,i}$ from there. Let $p_{n,i}'$ be the number of
$(n-i)$-multipartitions of $n$ such that every cell has an even number
of elements, e.g., $(2,2,4)$ is allowed, but not $(1,2,3,4)$.
\begin{corollary}
Let $B$ be an arbitrary filtered quantization of $\cO_{\CC^{2n}}^{D_n}$. Then, the
number of prime ideals of $B$ whose support has
codimension $2i$ in $V/D_n$ is at most 
\begin{equation}
p_{n,i}' + \sum_{j=0}^{i} d_j p_{n-j,i-j}.
\end{equation}
\end{corollary}

\subsection{The $\caD$-module $M(\CC^{2n}/D_n)$}\label{ss:typed-int-dmod}
Similarly to the case of symmetric powers of symplectic varieties in
\S \ref{ss:typea-int-dmod}, we may deduce the structure of $M(X)$ for
$X = \CC^{2n} / D_n$.

When $U$ is a vector space, let $\delta_{0 \in U}$ denote the
$\delta$-function $D_U$-module at the origin.  Let $q: \CC^{2n} \onto \CC^{2n}/D_n$ be the quotient map. 
Let $\Delta_i : \CC^2 \into (\CC^2)^{i}$ denote the diagonal
embedding.  Also, define the modified embedding $\Delta_i': \CC^2
\into (\CC^{2})^i$ by $\Delta_i'(x,y) = ((-x,-y),(x,y),\ldots,(x,y))$,
i.e., the composition of $-\Id \times \Id^{i-1}$ with $\Delta_i$.  
\begin{theorem}\label{dndmthm}
\begin{multline}
  M(X) \cong \\
  \bigoplus_{\underset{r_j \geq 1\, \forall j, 1 \leq i_1 < \cdots < i_k, r \geq
      0} {r + r_1 \cdot i_1 + \cdots + r_k \cdot i_k = n}}
  \HP_0(\cO_{\CC^{2r}}^{D_{r}}) \otimes q_* \bigl(\delta_{0 \in
    \CC^{2r}} \boxtimes
  \bigl((\Delta_{i_1})_*(\Omega_{\CC^2})^{\boxtimes r_1} \boxtimes
  \cdots \boxtimes (\Delta_{i_k})_*(\Omega_{\CC^2})^{\boxtimes
    r_k}\bigr)
  \bigr)^{\Stab} \\
  \oplus \!\!\! \bigoplus_{\underset{r_j \geq 1\, \forall j, 2 \leq i_1 < \cdots <
      i_k, 2 \mid i_j, \forall j}{r_1 \cdot i_1 + \cdots + r_k \cdot i_k = n}}
  q_*\bigl(
\bigl((\Delta_{i_1}')_*(\Omega_{\CC^2}) \boxtimes (\Delta_{i_1})_*(\Omega_{\CC^2})^{\boxtimes r_1-1}\bigr)
 \boxtimes  (\Delta_{i_2})_*(\Omega_{\CC^2})^{\boxtimes r_2} \boxtimes
  \cdots 
  \boxtimes (\Delta_{i_k})_*(\Omega_{\CC^2})^{\boxtimes
    r_k}\bigr)\bigr)^{\Stab}.
\end{multline}
\end{theorem}
Here, the superscript of $\Stab$ refers to the subgroup of $D_n$ which
preserves the support of the $\caD_{\CC^{2n}}$-module we are pushing
forward by $q$: for example, in the first big direct sum, this will be
the subgroup for each summand preserving the locus $\{(0,0)\}^{r}
\times \Delta_{i_1}(\CC^2)^{r_1} \times \cdots \times
\Delta_{i_k}(\CC^2)^{r_k}$. This group is explicitly
\[
D_n \cap \biggl(B_r \times \prod_{j=1}^k((\pm S_{i_j})^{r_j} \rtimes S_{r_j})\biggr),
\]
where here $\pm S_{i} \cong S_i \times \ZZ/2$ is the group generated by permutation matrices and $-\Id$.  One can express the second stabilizer in a similar way, and it is isomorphic to $\prod_{j=1}^k((\pm S_{i_j})^{r_j} \rtimes S_{r_j})$ (with the case $j=1$ of the product acting in a modified way so as to preserve the
locus $\Delta_{i_1}'(\CC^2) \times \Delta_{i_1}(\CC^2)^{r_1-1}$ rather than
$\Delta_{i_1}(\CC^2)^{r_1})$.

\subsection{Explicit computational results}
Using programs \cite{dnhp0-progs} written in Magma \cite{magma}, we
explicitly solved \eqref{e:xik-cond} for $n \leq 34$ (and
double-checked, for $n \leq 7$ and low enough degrees for $n \in \{8,
9\}$, that the result matches a direct computation of $\HP_0$ without
using Theorem \ref{dnhp0thm}).  The result is given in Figures \ref{f:expdn1} and \ref{f:expdn2}.

\begin{figure}
\begin{tabular}{c|c}
$n$ & $h(\HP_0(\cO_{\CC^{2n}}^{D_n});t^{\frac{1}{4}})$ \\[4 pt] \hline
\input{hs-dn-tabular1}
\end{tabular}
\caption{Poisson traces on type $D_n$ singularities
for $n \leq 19$\label{f:expdn1}}
\end{figure}

\begin{figure}
\begin{tabular}{c|c}
$n$ & $h(\HP_0(\cO_{\CC^{2n}}^{D_n});t^{\frac{1}{4}})$ \\[4 pt] \hline
\input{hs-dn-tabular2}
\end{tabular}
\caption{Poisson traces on type $D_n$ singularities
for $20 \leq n \leq 34$\label{f:expdn2}}
\end{figure}

\subsection{Proof of Theorem \ref{dnhp0thm}}
Set $A := \cO_{\CC^{2n}}^{D_n}$. (When we need $n$ to vary later on,
we will also denote $A$ by $A^{(n)}$.)  Fix a basis $x_1, \ldots, x_n,
y_1, \ldots, y_n$ of $(\CC^{2n})^*$ such that $\{x_i, y_j\} =
\delta_{ij}$ and $\{x_i, x_j\} = 0 = \{y_i, y_j\}$ for all $i, j$.
Decompose $A = A_+ \oplus A_-$ as the eigenspaces of the diagonal
matrix $T := \begin{pmatrix} -1 & & & \\ & 1 & & \\ & & \ddots & \\ &
  & & 1\end{pmatrix} \in \GL(\CC^n) \subseteq \Sp(\CC^{2n})$, in the
basis of the $x_i$ (or equivalently, any element of $B_n \setminus
D_n$).  Then, $A_+ = \cO_{\CC^{2n}}^{B_n}$ is the ring of polynomials
which are symmetric under the action of $S_n$ simultaneously on $x_i$
and $y_i$ and for which every monomial has an even sum of degrees in
the index-$i$ variables $x_i$ and $y_i$, for all $i$. Similarly, $A_-$
is the space of symmetric polynomials such that every monomial has
\emph{odd} total degree in $x_i$ and $y_i$, for all $i$. Note the
formula
\begin{equation} \label{am1fla}
A_- = \sum_{z_i \in \{x_i, y_i\}}
z_1 z_2 \cdots z_n A_+.
\end{equation}

We would like to compute $A/\{A,A\} = A_+/(\{A_+, A_+\} + \{A_-,
A_-\}) \oplus A_-/\{A_+, A_-\}$.
\begin{lemma} $A_- = \{A_+, A_-\}$. \label{am1vanl}
\end{lemma}
\begin{proof}Let $\symm(f) = \frac{1}{n!} \sum_{\sigma \in S_n}
  \sigma(f)$ be the symmetrization map.  We need to show that, for all
  monomials $x_1^{a_1} y_1^{b_1} \cdots x_n^{a_n} y_n^{b_n}$ such that
  $a_i + b_i$ is odd for all $i$, the symmetrization $\symm(x_1^{a_1}
  y_1^{b_1} \cdots x_n^{a_n} y_n^{b_n})$ is a sum of Poisson brackets.
  
  To do so, we consider a filtration on $A$ (which we will not label
  by integers) given by an ordering on monomials.  First, take the
  ordering on monomials in $\CC[x,y]$ of the form $x^a y^b > x^{a'}
  y^{b'}$ if either $a+b > a'+b'$ or $a+b=a'+b'$ and $a > a'$. Extend
  this to symmetrizations of
  monomials in $\CC[x_1, \ldots, x_n, y_1, \ldots, y_n]$:
  assuming that $x^{a_1} y^{b_1} \geq x^{a_2} y^{b_2} \geq \cdots \geq
  x_n^{a_n} y_n^{b_n}$ and similarly $x^{a_1'} y^{b_1'} \geq x^{a_2'}
  y^{b_2'} \geq \cdots \geq x_n^{a_n'} y_n^{b_n'}$, then we say that
  $\symm(x_1^{a_1} y_1^{b_1} \cdots x_n^{a_n} y_n^{b_n}) \geq \symm(x_1^{a_1'}
  y_1^{b_1'} \cdots x_n^{a_n'} y_n^{b_n'})$ if, for some $i$, $a_j=a_j'$
  and $b_j = b_j'$ for all $j < i$ and $x^{a_i} y^{b_i} > x^{a_i'}
  y^{b_i'}$.   

  For each degree $m \geq 0$, we will consider the induced filtration on the 
vector space $(A_-)_m$ in total degree $m$, given by the union
\begin{multline}
(A_-)_m = \bigcup_{\underset{a_i + b_i \text{ is odd}, \forall i}{a_1 + \cdots + a_n + b_1 + \cdots + b_n = m}} \Span\{\symm(x_1^{a_1'} y_1^{b_1'} \cdots x_n^{a_n'} y_n^{b_n'}) \text{ s.t. } \\
\symm(x_1^{a_1'} y_1^{b_1'} \cdots x_n^{a_n'} y_n^{b_n'}) > 
\symm(x_1^{a_1} y_1^{b_1} \cdots x_n^{a_n} y_n^{b_n})\},
\end{multline}
i.e., the map $(A_-)_m \to \gr (A_-)_m$ with respect to this
filtration takes the \emph{lowest symmetrized monomial} with respect
to the ordering on monomials, which has nonzero coefficient.

It suffices to show that $\gr \{A_-, A_+\}_m = \gr (A_-)_m$ for all $m
\geq 0$.  That is, for each monomial $x_1^{a_1} y_1^{b_1} \cdots
x_n^{a_n} y_n^{b_n}$ of total degree $m$ and with $a_i+b_i$ odd for
all $i$, we need to show that $\symm(x_1^{a_1} y_1^{b_1} \cdots
x_n^{a_n} y_n^{b_n}) + \text{higher terms} \in \{A_-, A_+\}_m$, for
some linear combination of greater symmetrizations of monomials of
total degree $m$.

So, assume $x^{a_i} y^{b_i} \geq x^{a_{i+1}} y^{b_{i+1}}$ for all $1 \leq i \leq n-1$.  We compute that
\begin{multline}
\{\symm(x_1^{a_1+1}y_1^{b_1}), \symm(y_1 \cdot x_2^{a_2} y_2^{b_2} \cdots x_n^{a_n} y_2^{b_n})\} \\ = \frac{1+c}{n}\cdot (a_1+1) \symm(x_1^{a_1} y_1^{b_1} \cdots x_n^{a_n} y_n^{b_n})
+ \text{higher terms},
\end{multline}
where $c = |\{i \in \{2,3,\ldots,n\} \mid (a_i, b_i) = (0,1)\}| \geq 0$.
\end{proof}
We conclude that
\begin{equation}
\HP_0(A) = \HP_0(\cO_{\CC^{2n}}^{B_n}) / \{A_-, A_-\}.
\end{equation}
Next, recall that, for every Poisson algebra $P$ which is Poisson
generated by elements $p_1, \ldots, p_k$, $\{P, P\} = \{\langle p_1,
\ldots, p_k \rangle, P\}$.  This is a result of the Jacobi identity
and the identity $\{ab, c\} + \{bc, a\} + \{ca, b\} = 0$.  Then, note
that $A_+ \subset A$ contains the copy of $\mathfrak{sl}_2$ spanned by
$\sum_i x_i^2$, $\sum_i y_i^2$, and $\sum_i x_i y_i$.  Here and below all sums over $i$ will range from $1$ to $n+1$ (not only $1$ to $n$) unless otherwise specified.

As a result of this and \eqref{am1fla}, $A$ is Poisson generated by
$A_+$ and the single element $y_1 y_2 \cdots y_n$.  Hence, we deduce
that
\begin{equation}
\HP_0(A) = \HP_0(\cO_{\hh \oplus \hh^*}^{B_n}) / \{y_1 y_2 \cdots y_n, A_-\}.
\end{equation}
As a result, there is a natural inclusion of dual spaces
$\HP_0(A)^* \subseteq \HP_0(\cO_{\CC^{2n}}^{B_n})^*$.  For convenience, when we allow $n$ to vary, let
$A^{(n)} := \cO_{\CC^{2n}}^{D_n} = A^{(n)}_+ \oplus
A^{(n)}_-$ be the above decomposition (note that $A^{(n)}_+ =
\cO_{\CC^{2n}}^{B_n}$).  In these terms, we showed that
$\HP_0(A^{(n)})^* \subseteq \HP_0(A_+^{(n)})^*$.
\begin{claim}
The image of the inclusion
\begin{equation}
\bigoplus_n \HP_0(A^{(n)})^* \subseteq \bigoplus_n \HP_0(A_+^{(n)})^*
\end{equation}
is a bigraded subalgebra.
\end{claim}
\begin{proof}
We have to show that, if $f \in \HP_0(A_+^{(m)})^* \subseteq (A_+^{(m)})^*$ and $g \in \HP_0(A_+^{(n)})^* \subseteq (A_+^{(n)})^*$ satisfy
\begin{equation}
f(\{y_1 \cdots y_m, A_-^{(m)}\}) = 0, \quad g(\{y_1 \cdots y_n, A_-^{(n)}\}) = 0,
\end{equation}
i.e., $f \in \HP_0(A^{(m)})^*$ and $g \in \HP_0(A^{(n)})^*$, then 
\begin{equation}
(f \cdot g) (\{y_1 \cdots y_{m+n}, A_-^{(m+n)}\}) = 0.
\end{equation}
This follows immediately from the Leibniz rule and the fact that, as
subspaces of $\cO_{\CC^{2(m+n)}} = T^{m+n}\cO_{\CC^2}$,
$A_-^{(m+n)} \subseteq A_-^{(m)} \otimes A_-^{(n)}$.
\end{proof}
% \begin{remark}
%   We did not really need Lemma \ref{am1vanl} to show that $\bigoplus_n
%   \HP_0(A^{(n)})^*$ is an algebra: the Leibniz rule similarly implies
%   that $\bigoplus_n (A^{(n)}_{-1} / \{A^{(n)}_1, A^{(n)}_{-1}\})^*$ is
%   a nonunital algebra (although Lemma \ref{am1vanl} shows that it is
%   zero).
% \end{remark}
We now explicitly describe the subalgebra $\bigoplus_n \HP_0(A^{(n)})^*
\subseteq \bigoplus_n \HP_0(A_+^{(n)})^*= \CC[s_i]_{i \geq 1}$.  This
depends on the choice of the $s_i$, each of which is canonical up to
scaling. We will make use of the construction of \cite{ESsym}, as we
recall in the proof.

Let us recall the definition of the functions $s_i$ from \cite[\S
4]{ESsym}. It is convenient to view $s_i$ as a degree-$i$ function on
$\CC[x^2,xy,y^2]$, i.e., $s_i(f) := s_i(f^{\otimes i})$. Since they are
homogeneous, the $s_i$ extend to continuous functions on the 
completion $\CC\dblsqbrs{x^2,xy,y^2}$.
It is proved
in \emph{op. cit.} that every $f$ in the completion $\CC\dblsqbrs{x^2,xy,
y^2}$ of $\cO_{\CC^2}^{\ZZ/2}$ with nonvanishing second derivative
(``nondegenerate'') is equivalent up to continuous Poisson
automorphisms (i.e., even symplectomorphisms of the formal disc) to a
unique element of the form
\begin{equation}
f \sim y^2 + s_1 +  s_2 x^2 +  s_3 x^4 + \cdots.
\end{equation}
Then, $s_i(f)$ is defined as the above coordinate in this normal
form. This extends uniquely to the entire pro-vector space $\CC\dblsqbrs{x^2,
xy, y^2}$ (no longer requiring the nonvanishing second derivative
condition) since the degenerate elements form a codimension-two
subspace.  These $s_i$ restrict to functions on $\CC[x^2,xy,y^2]$ and
have degree $i$ and weight $4-4i$, i.e., $s_i$ is a degree-$i$
polynomial and $s_i(f(\lambda x, \lambda y)) = \lambda^{4i-4}
s_i(f(x,y))$.

We need to consider the value of the functions $s_n$ on brackets of
the form $\{y_1 y_2 \cdots y_n, g\}$ for $g \in A^{(n)}_{-}$.  Since
$A^{(n)}_- = \Sym^n A^{(1)}_- = \langle x, y \rangle \CC[x^2,xy,y^2]$,
it follows that $A^{(n)}_-$ is spanned by elements of the form $g =
f^{\otimes n}$ for $f \in A^{(1)}_-$.  Then, we notice that
\begin{equation}
  \{y_1 y_2 \cdots y_n,
  f^{\otimes n}\} = \symm(n\{y, f\} \otimes (yf)^{\otimes (n-1)}).
\end{equation}
Thus, the subalgebra $\bigoplus_n \HP_0(A^{(n)})^* \subseteq \CC[s_i]_{i \geq 1}$ consists of those polynomials $F$ such that, working over $\CC[\varepsilon]/(\varepsilon^2)$,
\begin{equation}\label{Feqn1}
F(yf + \varepsilon \{y,f\}) - F(yf) = 0, \quad \forall f \in A_-^{(1)}.
\end{equation}
Next, we write $yf$ in normal form up to even formal
symplectomorphisms.  We claim that the result is of the form $(y +
h)(y - h)$, where $h \in x \CC\dblsqbrs{x^2}$.  Indeed, since $yf(0)=0$, the
result must lie in $y^2 + x^2 \CC\dblsqbrs{x^2}$, and since $\CC\dblsqbrs{x^2}$ admits
square roots in $\CC\dblsqbrs{x}$, the result must be of the form $y^2 - h^2$ for some
$h \in x\CC\dblsqbrs{x}$.  

The next step is to write $yf + \varepsilon \{y,f\}$ in normal form
up to even formal symplectomorphisms.  First let $\varphi$ be the
 aforementioned symplectomorphism satisfying
$\varphi(yf) = (y+h)(y-h)$.  Then, up to choice of $h$, $\varphi$ takes $y$ to $u(y + h)$
and $f$ to $u^{-1} (y - h)$, for some even unit $u \in \CC\dblsqbrs{x^2,xy,y^2}$.  Therefore,
\begin{multline}
\varphi(\{y, f\})
= \{u(y + h), u^{-1} (y - h)\} \\= 2\{h, y\} + u^{-1}(y + h)\{u, y -
h\} -u(y - h) \{u^{-1}, y + h\} \\ = 2\{h,y\} + u^{-1}(y + h)\{u, y - h\} + u^{-1}(y - h) \{u, y + h\}
\\ = 2 \{h, y\} + 2y \{\log(u), y\} -2h \{\log(u), h\} \\ = 2 \{h,y\} + \{\log(u), y^2 - h^2\}.
\end{multline}
Hence, $\varphi(yf + \varepsilon \{y,f\})= (y^2-h^2) + 2\varepsilon
\{h,y\} + \varepsilon \{\log(u), y^2 - h^2\}$.  Further, we may apply
the symplectomorphism $e^{-\varepsilon\ad(\log u)}$ and we obtain
$(y^2-h^2) + 2\varepsilon \{h,y\}$. Therefore, for $F \in \HP_0(A^{(n)})^*$,
\eqref{Feqn1} becomes
\begin{equation}\label{Feqn2}
  F(y^2-h^2 + 2 \varepsilon \{h,y\}) - F(y^2-h^2) = 0, \forall h \in x \CC\dblsqbrs{x^2}.
\end{equation}

Next, consider the $s_i$ as coordinate functions on the
infinite-dimensional affine space $y^2 + \CC\dblsqbrs{x^2}$, or just
$\CC\dblsqbrs{x^2}$ by deleting the $y^2$.  Then, the above equation says that
$F$ is annihilated by a particular (discontinuous) vector field up to
sign, which is supported on $x^2\CC\dblsqbrs{x^2}$, and at the point $g =
-h^2$ is the vector $2 h' = 2 (\sqrt{-g})'$.  For all $k \geq
0$, the square-root function is a regular
multivalued function on the locus $\CC^\times \cdot x^{2k} +
x^{2k+2}\CC\dblsqbrs{x^2}$, defined by 
\begin{equation}
\sum_{i \geq k} s_{i+1} x^{2i} \mapsto
\pm \sqrt{s_{k+1}} x^k Q(\frac{1}{s_{k+1}}\sum_{i \geq k+1}
s_{i+1}x^{2(i-k)}),
\end{equation}
 where $Q$ is the
Taylor series for $\sqrt{1+x}$. That is, it says that
\begin{equation}\label{skgreqn}
  \xi_k F|_{s_1 = s_2 = \cdots = s_{2k-1} = 0, s_{2k} \neq 0} = 0, \quad \forall k \geq 1,
\end{equation}
with $\xi_k$ as in the statement of the theorem.

\subsection{Proof of Theorem \ref{dndmthm}}\label{ss:dndmthm-proof}
Let $X := \CC^{2n}/D_n$. By \cite[Theorem 3.1]{ESdm} and its proof, it
suffices to analyze $M(X)$ in the formal neighborhood of each
symplectic leaf. These leaves are the image under $\CC^{2n} \onto X$ of
two types of partitions:
\begin{enumerate}
\item[(i)] Partitions $\{1,\ldots,n\} := I \sqcup J_1 \sqcup \cdots \sqcup
J_{\ell}$, with the leaf an open subset of the locus where $x_i = 0 =
y_i$ for all $i \in I$, and $x_i = x_j$ and $y_i = y_j$ for all $i,j
\in I_k$, for all $1 \leq k \leq \ell$;
\item[(ii)] Partitions $\{1,\ldots,n\} := J_1 \sqcup \cdots \sqcup
J_{\ell}$ with each $|J_i|$ even and $1,2 \in J_1$; the leaf is 
an open subset of the locus where $x_i = x_j$ and $y_i = y_j$ for all $i,j
\in I_k$, for all $1 \leq k \leq \ell$, except that $x_1 = -x_2$ and $y_1 = -y_2$.
\end{enumerate}
Namely, the leaves are the complement in such loci of properly
contained such loci, i.e., those corresponding to partitions obtained
by joining some cells of the given partition (and if $I$ joins with
any other cell $J_i$, the new label must remain $I$).  The computation
therefore reduces to the cases $\{1,\ldots,n\} = I$ or $\{1,\ldots,n\}
= J_1$.  In the former case, the local system is just a multiple of
the delta-function local system at zero, whose multiplicity must be
$\dim \HP_0(\cO_X)$, and in the latter case, the problem reduces to the
computation of Theorem \ref{sympdmthm}.

\section{Direct proof of Theorem \ref{snp1thm}, by T. Schedler} \label{s:direct-proof-snp1thm}
%Let $\hh$ be the reflection representation of $S_n$, so $\dim\ \hh =
%n-1$.  Let $V := \hh \oplus \hh^*$.  
We need to show that $\HP_0(\cO_V^{S_n}, \cO_V^{S_n}) = \CC$.  In
fact, we will show the stronger result that $\HP_0(\cO_V^{S_{n}},
\cO_V^{S_{n-1}}) = \CC$ (note that this is also what Lemma
\ref{l:hidist} proves). It will be convenient to remember
that $\{\cO_V^{S_n}, \cO_V\}^{S_{n-1}} = \{\cO_V^{S_n},
\cO_V^{S_{n-1}}\}$.
% It is clear that $\HP_0(\CC[V]^{S_n}, \CC[V]^{S_{n-1}})$ is graded
% by the usual grading $|V^*|=1$ on $\CC[V]$, and the Poisson bracket
% has degree $-2$.

In degree zero, we clearly get $\HP_0(\cO_V^{S_n},\cO_V^{S_{n-1}})_0 =
\CC$.  So it suffices to show that the positively-graded part
vanishes, i.e., $\HP_0(\cO_V^{S_n}, \cO_V^{S_{n-1}})_{> 0} = 0$.

It will be helpful to explicitly write $V$ in terms of
coordinates. Let $V = \mathfrak{h} \oplus \mathfrak{h}^*$ where
$\mathfrak{h} \cong \CC^{n-1}$ is the reflection representation. We
can consider $\mathfrak{h} \subseteq \CC^n$ to be the subset where all
coordinates sum to zero.  Hence we can write $\cO_V = \CC[x_1, \ldots,
x_n, y_1, \ldots, y_n] / (x_1 + \cdots + x_n, y_1 + \cdots + y_n)$.
Moreover,  we can choose coordinates so that the permutation action of
$S_n$ on $x_1, \ldots, x_n$ and on $y_1, \ldots, y_n$ is by the usual
action (by simultaneous permutations of indices using the same
permutation), and the Poisson bracket is given by
\[
\{x_i, y_j\} = \delta_{ij} - \frac{1}{n}.
\]
(Note that the $-\frac{1}{n}$ is required here because, for instance,
$x_1 + \cdots + x_n = 0$.)

Consider the sub-Lie algebra of $\cO_V^{S_n}$ spanned by $\sum_{i=1}^n
x_i^2$, $\sum_{i=1}^n y_i^2$, and $h := \sum_{i=1}^n x_i y_i$.  This
is isomorphic to $\mathfrak{sl}_2$, and we will simply call it
$\mathfrak{sl}_2$. Moreover, the action of $S_n$ commutes with the
action of $\mathfrak{sl}_2$.

Since the adjoint action of $\mathfrak{sl}_2$ preserves degree on
$\cO_V$, in each degree we obtain a semisimple representation. Hence,
$\{\cO_V^{S_n}, \cO_V\}$ contains the sum of all nontrivial
$\mathfrak{sl}_2$-representations in $\cO_V$, $\{\mathfrak{sl}_2,
\cO_V\}$.  Because the action of $S_n$ commutes with that of
$\mathfrak{sl}_2$, also $\{\cO_V^{S_n}, \cO_V^{S_{n-1}}\}$ contains
the $S_{n-1}$-invariants of the sum of nontrivial
$\mathfrak{sl}_2$-representations, $\{\mathfrak{sl}_2, \cO_V\}^{S_{n-1}}$.

Next, for any two finite-dimensional $\mathfrak{sl}_2$-representations
$W$ and $W'$, if $w' \in W'$ is a highest (or lowest) weight vector
for $h \in \mathfrak{sl}_2$, it is easy to see that $W \otimes w'$
generates $W \otimes W'$ as a $\mathfrak{sl}_2$-representation (e.g.,
one can assume $W$ is irreducible, and then show that all tensor
products $w_1 \otimes w_2$ of $h$-weight vectors are generated, by
induction on the weights). Since $y_n^k$ is a highest weight vector
for the representation $\CC[y_n]_k$ (the subscript denotes degree
$k$),
% Since $\CC[x_i, y_i]$ is a $\mathfrak{sl}_2$-subrepresentation, and the trivial representation occurs only once (spanned by $1 \in \CC[x_i,y_i]$),
it follows that, for all $1 \leq j < n$,
\begin{equation}\label{e:sl2-trick}
  (\CC[x_1,\ldots,x_j,y_1,\ldots,y_j] \cdot \CC[y_n]) + 
  \{\mathfrak{sl}_2, \cO_V\} \supseteq 
  \CC[x_1,\ldots,x_j,y_1,\ldots,y_j] \cdot \CC[x_n,y_n].
\end{equation}

Let $\symm: \cO_V \rightarrow \cO_V^{S_{n}}$ be the symmetrization
map, and $\symm': \cO_V \rightarrow \cO_V^{S_{n-1}}$ be the
symmetrization for the subgroup $S_{n-1} \subset S_n$.

In view of the above, it suffices to prove that, for all $1 \leq j < n$
and all $a_1, \ldots, a_{j}, b_1, \ldots, b_{j} \geq 0$ and $m \geq 0$ such
that the $a_i, b_i$, and $m$ are not all zero,
\begin{equation} \label{indcl}
\symm'\bigl(x_1^{a_1} y_1^{b_1} \cdots x_{j}^{a_{j}} y_{j}^{b_{j}} \cdot
y_n^m \bigr) \in \{\cO_V^{S_n}, \cO_V^{S_{n-1}}\}.
\end{equation} 
It suffices to fix the total degree $d = \sum_i (a_i + b_i) + m \geq
1$, which we fix from now on. Note that, once we prove \eqref{indcl}
for all $a_i, b_i, m$ such that the total degree is $d$, it follows
from \eqref{e:sl2-trick} that
\begin{equation}\label{indcl2}
\symm'\CC[x_1,y_1,\ldots,x_j,y_j,x_n,y_n]_d \subseteq 
\{\cO_V^{S_n}, \cO_V^{S_{n-1}}\},
\end{equation}
for the same value of $j$ (and the same degree $d$).  Here the
subscript of $d$ denotes the part of total degree $d$. Now, we will
prove \eqref{indcl}, for fixed total degree $d$, by a double
induction: we induct on $j$, and for each value of $j$, we also
perform a reverse induction on $m$ (beginning with $m=d$, the total
degree).

For the base case(s) of the (double) induction, for any value of $j$,
it is enough to show that $\symm'(y_n^d)$ is in $\{\cO_V^{S_n},
\cO_V^{S_{n-1}}\}$.  This follows because $y_n^d$ generates a
nontrivial representation of $\mathfrak{sl}_2$ in $\cO_V$, since $d
\geq 1$.

The inductive step follows from the computation
\begin{multline}
\symm'(x_1^{a_1}y_1^{b_1} \cdots x_{j}^{a_{j}}y_j^{b_j} y_{n}^{a_{n}}) \\= -\frac{n^2}{(n-j)(a_1+1)(a_{n}+1)}\symm'\Bigl(\{\symm(x_1^{a_1+1}y_1^{b_1}),
x_2^{a_2} y_2^{b_2} \cdots x_{j}^{a_{j}}y_j^{b_j} y_{n}^{a_{n}+1}\} + \text{h.o.t.}\Bigr),
\end{multline}
where ``h.o.t.''=''higher order terms'' refers to a linear combination
of monomials  with fewer
indices appearing (in this case, the only variables which occur in this
part of the sum will be $x_2, \ldots, x_{j},
y_2, \ldots, y_j, x_n$, and $y_{n}$), or where the exponent of $y_n$ appearing is greater (in this case, it will be $y_n^{a_n+1}$).  These are already in
$\{\cO_V^{S_n}, \cO_V^{S_{n-1}}\}$ by hypothesis.
This completes the induction.
\bibliographystyle{amsalpha}
\bibliography{master}
\end{document}